\documentclass[twoside,leqno,10pt]{article}
\usepackage{graphics,ifthen}
\usepackage{latexsym}
\usepackage{amsmath}
\usepackage{amssymb}
\usepackage{mathrsfs}

\begin{document}
\input{latexP.sty}
\input{referencesP.sty}
\input epsf.sty

\def\ind{\stackrel{\mathrm{ind}}{\sim}}
\def\iid{\stackrel{\mathrm{iid}}{\sim}}

\def\Definition{\stepcounter{definitionN}\
    \Demo{Definition\hskip\smallindent\thedefinitionN}}
\def\EndDefinition{\EndDemo}
\def\Example#1{\Demo{Example [{\rm #1}]}}
\def\EndExample{\qed\EndDemo}
\def\Category#1{\centerline{\Heading #1}\rm}
\
%% Paper specific definitions
\def\e{\text{\hskip1.5pt e}}
\newcommand{\eps}{\epsilon}
\newcommand{\proof}{\noindent {\bf Proof:\ }}
\newcommand{\remarks}{\noindent {\bf Remarks:\ }}
\newcommand{\note}{\noindent {\bf Note:\ }}
\newcommand{\examp}{\noindent {\bf Example:\ }}
\newcommand{\Lower}[2]{\smash{\lower #1 \hbox{#2}}}
\newcommand{\ben}{\begin{enumerate}}
\newcommand{\een}{\end{enumerate}}
\newcommand{\bi}{\begin{itemize}}
\newcommand{\ei}{\end{itemize}}
\newcommand{\hp}{\hspace{.2in}}

\newtheorem{lw}{Proposition 3.1, Lo and Weng (1989)}
\newtheorem{thm}{Theorem}[section]
\newtheorem{defin}{Definition}[section]
\newtheorem{prop}{Proposition}[section]
\newtheorem{lem}{Lemma}[section]
\newtheorem{cor}{Corollary}[section]
\newcommand{\rb}[1]{\raisebox{1.5ex}[0pt]{#1}}
\newcommand{\mc}{\multicolumn}
%Mathrsfs Font
\newcommand{\Bcr}{\mathscr{B}}
\newcommand{\Ucr}{\mathscr{U}}
\newcommand{\Gcr}{\mathscr{G}}
\newcommand{\Dcr}{\mathscr{D}}
\newcommand{\CS}{\mathscr{C}}
\newcommand{\Fcr}{\mathscr{F}}
\newcommand{\Icr}{\mathscr{I}}
\newcommand{\Lcr}{\mathscr{L}}
\newcommand{\Mcr}{\mathscr{M}}
\newcommand{\Ncr}{\mathscr{N}}
\newcommand{\Pcr}{\mathscr{P}}
\newcommand{\Qcr}{\mathscr{Q}}
\newcommand{\Scr}{\mathscr{S}}
\newcommand{\Tcr}{\mathscr{T}}
\newcommand{\Xcr}{\mathscr{X}}
\newcommand{\Vcr}{\mathscr{V}}
\newcommand{\Ycr}{\mathscr{Y}}
\newcommand{\qcr}{\mathscr{q}}
%Mathbb Font
\newcommand{\E}{\mathbb{E}}
\newcommand{\F}{\mathbb{F}}
\newcommand{\I}{\mathbb{I}}
\newcommand{\Q}{\mathbb{Q}}
\newcommand{\X}{\mathbb{X}}
\newcommand{\Pe}{\mathbb{P}}
\newcommand{\M}{\mathbb{M}}
\newcommand{\Wbb}{\mathbb{W}}

\def\Beta{\text{Beta}}
\def\Dir{\text{Dirichlet}}
\def\DP{\text{DP}}
\def\P{{\bf p}}
\def\fhat{\widehat{f}}
\def\GA{\text{gamma}}
\def\ind{\stackrel{\mathrm{ind}}{\sim}}
\def\iid{\stackrel{\mathrm{iid}}{\sim}}
\def\J{{\bf J}}
\def\K{{\bf K}}
\def\min{\text{min}}
\def\N{\text{N}}
\def\p{{\bf p}}
\def\U{{\bf U}}
\def\W{{\bf W}}
\def\S{{\bf S}}
\def\T{{\bf T}}
\def\y{{\bf y}}
\def\t{{\bf t}}
\def\m{{\bf m}}
\def\X{{\bf X}}
\def\Y{{\bf Y}}
\def\tps{\mbox{\scriptsize ${\theta H}$}}   %   smaller "\psi"-vector
\def\ups{\mbox{\scriptsize ${P_{\theta}(g)}$}}   %   smaller "\psi"-vector
\def\vps{\mbox{\scriptsize ${\theta}$}}   %   smaller "\psi"-vector
\def\vups{\mbox{\scriptsize ${\theta >0}$}}   %   smaller "\psi"-vector
\def\hps{\mbox{\scriptsize ${H}$}}   %   smaller "\psi"-vector
\def\rps{\mbox{\scriptsize ${(\theta+1/2,\theta+1/2)}$}}   %   smaller "\psi"-vector
\def\sps{\mbox{\scriptsize ${(1/2,1/2)}$}}   %   smaller "\psi"-vector

\newcommand{\reals}{{\rm I\!R}}
\newcommand{\PR}{{\rm I\!P}}
\def\Z{{\bf Z}}
\def\yy{{\mathcal Y}}
\def\rr{{\mathcal R}}
\def\BP{\text{beta}}
\def\ts{\tilde{t}}
\def\js{\tilde{J}}
\def\gs{\tilde{g}}
\def\fs{\tilde{f}}
\def\ys{\tilde{Y}}
\def\ps{\tilde{\mathcal {P}}}

\def\Report{Lancelot F. James}
\def\Author{OU-Gamma}
\pagestyle{myheadings}
\markboth{\Author}{\Report}
\thispagestyle{empty}

\bct\Heading  Laws and Likelihoods for Ornstein Uhlenbeck-Gamma
and other BNS OU Stochastic Volatilty models with extensions.
\lbk\lbk\smc Lancelot F. James\footnote{\eightit \rm Supported in
part by grants HIA05/06.BM03 and DAG04/05.BM56 of the HKSAR.
\eightit Corresponding authors address.
                \rm The Hong Kong University of Science and Technology,
Department of Information and Systems Management, Clear Water Bay,
Kowloon, Hong Kong.
\rm lancelot\at ust.hk\\
\indent\eightit Keywords and phrases.
                \rm
          Bessel Functions,
          Dilogarithm function,
          Dirichlet Process,
          Ornstein-Uhlenbeck Process,
          Perfect Sampling,
          Stochastic Volatility,
          Weber-Sonine Formula.
          }
\lbk\lbk \BigSlant The Hong Kong University of Science and
Technology\rm \lbk %(\today)%
\ect \Quote In recent years there have been many proposals as
flexible alternatives to  Gaussian based continuous time
stochastic volatility models. A great deal of these models employ
positive L\'evy processes. Among these are the attractive
non-Gaussian positive Ornstein-Uhlenbeck~(OU) processes proposed
by Barndorff-Nielsen and Shephard~(BNS) in a series of papers. One
current problem of these approaches is the unavailability of a
tractable likelihood based statistical analysis for the returns of
financial assets. This paper, while focusing on the BNS models,
develops general theory for the implementation of statistical
inference for a host of models. Specifically we show how to reduce
the infinite-dimensional process based models to finite, albeit
high, dimensional ones. Inference can then be based on Monte Carlo
methods. As highlights, specific to BNS we show that an OU process
driven by an infinite activity Gamma process, that is an
OU-$\Gamma$, exhibits unique features which allows one to exactly
sample from relevant joint distributions. We show that this is a
consequence of the OU structure and the unique calculus of Gamma
and Dirichlet processes. Owing to another connection between
Gamma/Dirichlet processes and the theory of Generalized Gamma
Convolutions~(GGC) we identify a large class of models, we
call~(FGGC), where one can perfectly sample marginal distributions
relevant to option pricing and Monte Carlo likelihood analysis.
This involves a curious result,  we establish as Theorem 6.1. We
also discuss analytic techniques and candidate densities for
Monte-Carlo procedures which can be applied to more general
classes of models. \EndQuote
%\baselineskip14pt
%\begin{document}
\rm
%\newpage
\tableofcontents
\section{Introduction}
Barndorff-Nielsen and Shephard~(2001a, b)(BNS) introduce a class
of continuous time stochastic volatility~(SV) models that allows
for more flexibility over Gaussian based models such as the
Black-Scholes model[see Black and Scholes~(1973) and
Merton~(1973)]. Their proposed SV model is based on the following
differential equation, \Eq dx^{*}(t)=(\mu+\beta
v(t))dt+v^{1/2}(t)dw(t) \label{BNS}\EndEq where $x^*(t)$ denotes
the log-price level, $w(t)$ is Brownian motion, and independent of
$w(t)$, $v(t)$ is a stationary Non-Gaussian Ornstein-Uhlenbeck
(OU) process which models the {\it instantaneous volatility}. This
latter point is equivalent to the fact that for $\lambda>0$,
$$
v(t)={\mbox e}^{-\lambda t}v(0)+{\mbox e}^{-\lambda
t}\int_{0}^{t}{\mbox e}^{\lambda y}Z(d\lambda y)
$$
and arises as the solution of the following differential equation,
$$
dv(t)=-\lambda v(t)+dZ(\lambda t).
$$
In the above framework $Z$ is a positive homogeneous process,
otherwise known as a \emph{subordinator}, on $[0,\infty)$ and
$v(0)$ is an arbitrary positive random variable independent of
$Z.$ That is $Z(t):=\int_{0}^{t}Z(dy)$ is a stationary process,
with $Z(0)=0,$ and its distribution specified by its Laplace
transform for each $\omega>0$, \Eq \E[{\mbox e}^{-\omega
Z(t)}]={\mbox e}^{-t\psi(\omega)} \label{Laplace}\EndEq where
$\psi(\omega)=\int_{0}^{\infty}(1-{\mbox e}^{-s\omega})\rho(ds)$,
is often called the \emph{L\'evy exponent} of an infinite
divisible random variable equivalent in distribution to $Z(1)$,
and $\rho$ is its corresponding L\'evy density. Either of these
characterizes the distribution of the process $Z.$ Importantly, it
is obvious from~\mref{Laplace}, that one does not need explicit
knowledge of $\rho$ to calculate $\psi.$ Note further that if we
wish $v(t)$ to be {\it stationary} it is necessary to choose
$v(0)\overset {d}=\int_{-\infty}^{0}{\mbox e}^{s}Z^{*}(ds)$, where
$Z^{*}$ is independent of $Z$ but otherwise has the same law.

The  model described above is an extension of the Black-Scholes or
Samuelson model which arises by replacing $v$ with a fixed
variance, say $\sigma^{2}$. The additional innovation in BNS is
that modeling volatility as a random process, $v(t)$, rather than
a random variable, not only allows for heavy-tailed models, but
additionally induces serial dependence. This serial dependence is
used to account for a clustering affect referred to as {\it
volatility persistence}. The work of Carr, Geman, Madan, and
Yor~(2003) discuss this point further. See also Duan~(1995) and
Engle~(1982) for different approaches to this type of phenomenon.
The model of BNS has gained a great deal of interest with some
related works including Carr, Geman, Madan, and Yor~(2003),
Barndorff-Nielsen and Shephard~(2003), Eberlein~(2001), Nicolato
and Venardos~(2001), Benth, Karlsen, and Reikvam~(2003). See also
the discussion section in~Barndorff-Nielsen and Shephard~(2001a).
See Carr and Wu~(2004) and Duffie, Pan and Singleton~(2000) for
many other models.

Note that the log price at time $t$ is $ x^{*}(t)=\mu t+\beta
\tau(t)+\tau^{1/2}(t)w(t) $ where $$
\tau(t)=\int_{0}^{t}v(s)ds=\lambda^{-1}[(1-{\mbox e}^{-\lambda
t})v(0)+\int_{0}^{t}(1-{\mbox e}^{-\lambda(t-y)})Z(d\lambda y)]$$
is referred to as a \emph{integrated} OU process and models the
integrated variance. Quantities of interest are often based on the
aggregate returns, for $s<t$, $x^{*}(t)-x^{*}(s)$ which involves
 \Eq
\tau(t)-\tau(s)=\lambda^{-1}[(1-{\mbox
e}^{-\lambda(t-s)})v(s)+\int_{s}^{t}(1-{\mbox
e}^{-\lambda(t-y)})Z(d\lambda y)] \label{tdiff}\EndEq where again
importantly, $ v(s)=({\mbox e}^{-\lambda s}v(0)+\int_{0}^{s}{\mbox
e}^{-\lambda(s-y)}Z(d\lambda y)).$

Barndorff Nielsen and Shephard~(2001a, Section 5.4.1 and 6.2) show
that laws related to the random functions \Eq (Z(\lambda t),
{\mbox e}^{-\lambda t}\int_{0}^{t}{\mbox e}^{\lambda y}Z(d\lambda
y)) \label{pair1} \EndEq play a key role both in option pricing
and likelihood estimation. Specifically option pricing requires
some type of description of the distribution of
\Eq\int_{s}^{t}(1-{\mbox e}^{-\lambda(t-y)})Z(d\lambda
y)\overset{d}=\int_{0}^{\Delta}(1-{\mbox
e}^{-y})Z(dy)\label{Option},\EndEq for $\Delta=(t-s)>0.$ Although
the density of \mref{Option} is not often known in a nice closed
form one can apply inversion techniques via its characteristic
function or Laplace transform which is described in BNS (2001a,
2003).

However, as seen in BNS (2001a, 5.4) it is a rather challenging
problem to find tractable approaches to statistical analysis of
likelihood models based on $n$ aggregate returns, $X_{i}=
x^{*}({i\Delta})-x^{*}((i-1)\Delta)$ over periods of time
$[(i-1)\Delta,i\Delta]$ for $i=1,\ldots, n$ and $\Delta>0$.[This
framework can be extended to intervals of varying lengths say
$\Delta_{i}$]. These models are based on the unobserved {\it
actual variances} $\tau_{i}=\tau(i\Delta)-\tau((i-1)\Delta)$ for
$i=1,\ldots, n.$ It is easy to see that one may write \Eq
X_{i}=\mu\Delta+\beta \tau_{i}+\tau^{1/2}_{i}\epsilon_{i}
\label{BNSx}\EndEq where $\epsilon_{i}$ for $i=1,\ldots n$ are
independent standard Normal random variables. Hence it follows
that conditional on $(\tau_{1},\ldots, \tau_{n})$ the $X_{i}$ are
independent Normal random variables with unknown mean $\mu\Delta
+\beta \tau_{i}$ and variance $\tau_{i}.$ The major obstacle to
tractable analysis of such models is that in general the joint
distribution of $(\tau_{1},\ldots,\tau_{n})$ is rather complex.
BNS~(2001a, 5.4) show that statistical inference can be done if
one were able to sample, \emph{efficiently},  $n$ iid copies of
the pair \Eq (Z(\lambda \Delta),{\mbox e}^{-\lambda
\Delta}\int_{0}^{\Delta}{\mbox e}^{\lambda y}Z(d\lambda y)).
\label{pair2}\EndEq The problem is that is not obvious how to deal
with the joint distributional behavior of the above
pair~\mref{pair2}. This is in contrast to the option pricing
problem which essentially involves the distribution of a single
random variable. A generic, theoretically all purpose,  approach
is to use an infinite series representation. Several MCMC
procedures, based on variations of this idea, have been proposed
to handle subclasses of these models requiring simulation of
points from random processes. See for instance, Roberts,
Papaspiliopoulos and Dellaportas~(2004) and Griffin and
Steel~(2005), who use compound Poisson process specifications for
$Z,$ and the discussion section in Barndorff-Nielsen and
Shephard~(2001a). For approaches to other types of models see for
instance Eraker, Johannes, and Polson~(2003).

While these methods have their attractive points they do not
provide exact solutions for cases where $Z$ is an \emph{infinite
activity} process, such as a Gamma process or more generally a
Generalized Inverse Gamma (GIG) process. Moreover these methods
are computationally non-trivial and further work needs to be done
to assess their accuracy for different processes. Another
important point is that they cannot be used if one does not have
specific knowledge of the L\'evy density associated with $Z.$ This
excludes for instance the case where $Z$ is based on a Pareto or
LogNormal distribution.
\subsection{Proposal and outline}
This paper focuses on several subtopics related to the issues
above. In particular we discuss methods that avoid working
directly with infinite dimensional components. First, perhaps most
remarkably, we will show that if one chooses $Z$ to be a Gamma
process then one can sample exactly random variables based on the
pair in~\mref{pair2} and~\mref{pair1}. In addition, we will be
able to derive the explicit density of certain quantities which is
also relevant to option pricing. Curiously we will show that the
explicit densities depend on the {\it dilogarithm} function $$
Li_{2}(x):=-\int_{0}^{x}\frac{\log(1-u)}{u}du:=\sum_{k=1}^{\infty}\frac{x^{k}}{k^{2}}$$
The dilogarithm function is a well-studied special function that
arises often in a variety of contexts. See for instance
Maximon~(2003) and Flajolet and Sedgewick~(2006). This leads to an
explicit description of the relevant $\tau_{i}$ in terms of sums
of independent random variables which allows one to perform
likelihood estimation based on sampling $2n$ iid random random
variables as well as the independent random variable $v(0)$. We
then easily extend this framework to possibly random observation
times. An important point is that these results allow one to also
use $\tau$ in other likelihood models not discussed in BNS~(2001a,
5.4). These facts have not been pointed out in the literature.
They are derived from the unique properties of the Gamma/Dirichlet
process calculus wherein we are able to exploit a, not immediately
obvious, connection to Dirichlet Process mean functionals. In as
much, the seminal work of Cifarelli and Regazzini~(1990) and the
perfect simulation methods discussed in Guglielmi, Holmes and
Walker~(2002) play a key role. The corresponding OU process $v(t)$
is known as OU-$\Gamma$ process. This should not be confused with
the often discussed $\Gamma$-OU process where the BDLP is a
compound Poisson process and $v(t)$ has Gamma distributed
marginals. Also, our results suggest that one could simply use the
OU-$\Gamma$ as a building block for more intricate models.

Secondly all the properties that we exploit for the Gamma case do
not extend to other OU models. However we show that the ability to
perfectly sample the marginal distributions of quantities relevant
to option pricing and likelihood estimation extends to a large
class of models where $Z$ is a Generalized Gamma Convolution(GGC).
We call these models finite GGC or (FGGC). A highlight of this
paper related to this class of models is Theorem 6.1.

Although we shall focus primarily on the BNS OU models,  we note
that  there are many others which can be found for instance in
Carr and Wu~(2004)[see also Carr, Geman, Madan and Yor~(2003)]. As
such we shall employ an analytical technique which leads to an
expression of the relevant likelihoods in terms of an
$n$-dimensional Fourier-Cosine integral. This technique is loosely
based on the ideas in James~(2005b). Multidimensional
Fourier-Cosine integral appear often in various fields including
physics. We then focus on ingredients necessary to carry out Monte
Carlo procedures which are known to be well suited to
approximating high-dimensional integrals. More details may be
obtained from the provided table of contents.

\Remark Throughout, when appropriate,  we will be describing the
law of a generic positive random variable $W$ by its corresponding
L\'evy exponent defined as
$$
-\log E[{\mbox e}^{-\omega W}]
$$
We will use the notation $\Delta$ as an arbitrary positive
distance between two points. We shall specify its value when
necessary, i.e. $\Delta=t-s$, $\Delta=t$ and so on. We will also
often use the notation $a= \lambda \Delta.$\EndRemark
\section{Preliminaries}
This paper utilizes results from several linked but not often
jointly studied areas. We anticipate that the average reader will
be familiar with some but not all of the topics. As such we
provide some details that we shall exploit. The majority of the
discussion in sections 2.1-2.2 may be found in BNS~(2001a,b,
2003). Section 2.3 is again a blend of ideas from several fields.
\subsection{Some more preliminary OU results}
We will describe the distribution of pertinent quantities via
their L\'evy exponents, and discuss the basic structure of the
likelihood. First note that for any positive $g$ on $[0,\infty)$,
we may define a random variable
$Z(g):=\int_{0}^{\infty}g(x)Z(dx)$. Moreover it is fairly
well-known that the L\'evy exponent of $Z(g)$ is given by
$\int_{0}^{\infty}\psi(\omega g(x))dx.$  It is clear that all the
OU related processes that we encounter are representable as some
$Z(g)$ where $g(x)$ is readily identified. Using this fact or
consulting directly BNS(2001a,b, 2003) one has that the L\'evy
exponent of the quantity in~\mref{Option} is \Eq \int_{{\mbox
e}^{-\lambda \Delta}}^{1}\psi(\omega(1-u))u^{-1}du
\label{OLevy}\EndEq for $\Delta=(t-s)>0$. The L\'evy exponent of
${\mbox e}^{-\lambda t}\int_{0}^{t}{\mbox e}^{\lambda y}Z(d\lambda
y)\overset{d}={\mbox e}^{-\lambda \Delta}\int_{0}^{\Delta}{\mbox
e}^{\lambda y}Z(d\lambda y)$ for $\Delta =t$ is,$$\int_{{\mbox
e}^{-\lambda \Delta}}^{1}\psi(\omega u)u^{-1}du $$ If we wish to
choose $v(t)$ stationary then the L\'evy exponent of $v(0)$ must
be \Eq \int_{0}^{\infty}\psi(\omega {\mbox
e}^{-s})ds=\int_{0}^{1}\psi(\omega u)u^{-1}du.\label{VLevy}\EndEq
\subsection{BNS Likelihood model}
The model of Barndorff-Nielsen and Shephard~(2001a, section 5.4)
translates into a likelihood based model as follows.
 Let ${X_{i}}$ for $i=1,\ldots, n$ denote a sequence of aggregate
returns of the log price of a stock observed over intervals of
length $\Delta>0$, described in ~\mref{BNSx}. Suppose additionally
the $Z$ depends on unknown parameters $\varrho$. The likelihood of
the model depends on unknown parameters
$\vartheta=(\mu,\beta,\lambda,\varrho)$ and as stated before the
$X_{i}|\vartheta,\tau$ are iid Normal random variables. Ideally
one is interested in estimating $\vartheta$ based on the
likelihood \Eq \Lcr(\X|\vartheta)=\int_{{\mathbb R}^{n}_{+}}
\[\prod_{i=1}^{n}\phi(X_{i}|\mu\Delta+\beta
\tau_{i},\tau_{i})\]f(\tau_{1},\ldots,\tau_{n}|\varrho,\lambda)d\tau_{1},\ldots,d\tau_{n}
\label{BNSlik} \EndEq where, setting $A_{i}=(X_{i}-\mu\Delta)$,
and ${\bar A}=n^{-1}\sum_{i=1}^{n}A_{i}$,
$$\phi(X_{i}|\mu\Delta+\beta
\tau_{i},\tau_{i})={\mbox e}^{A_{i}\beta
}\frac{1}{\sqrt{2\pi}}\tau^{-1/2}_{i}{\mbox
e}^{-A^{2}_{i}/(2\tau_{i})}{\mbox e}^{-\tau_{i}\beta^{2}/2}
$$
denotes a Normal density. The quantity
$f(\tau_{1},\ldots,\tau_{n}|\varrho,\lambda)$ denotes the joint
density of the integrated volatility based on the intervals
$[(i-1)\Delta,i\Delta]$ for $i=1,\ldots, n$. Barndorff-Nielsen and
Shephard~(2001a) note that the likelihood is intractable and hence
makes exact inference difficult. The apparent intractability is
attributed to the complex nature of
$f(\tau_{1},\ldots,\tau_{n}|\varrho,\lambda)$ which is derived
from a random measure. Specifically, the BNS models complexities
arises from the following structure of the $\tau_{i}$.
From~\mref{tdiff} one has for the BNS model \Eq
\lambda\tau_{i}=(1-{\mbox e}^{-\lambda
\Delta})v((i-1)\Delta)+\int_{(i-1)\Delta}^{i\Delta}(1-{\mbox
e}^{-\lambda(i\Delta-y)})Z(d\lambda y)\label{tausimp1} \EndEq
where importantly for $r_{j}={\mbox e}^{\lambda j\Delta}$, and $$
\textsc{O}_{j}=\int_{(j-1)\Delta}^{j\Delta}{\mbox
e}^{-\lambda(j\Delta-y)})Z(d\lambda y),$$ $ v((i-1)\Delta)={\mbox
e}^{-\lambda (i-1)\Delta}[v(0)+\sum_{j=1}^{i-1}r_{j}
\textsc{O}_{j}]. $ It is not difficult to see that the $O_{j}$ are
iid for $j=1,\ldots, n$ but are correlated with corresponding
terms $$ \int_{(j-1)\Delta}^{j\Delta}(1-{\mbox
e}^{-\lambda(j\Delta-y)})Z(d\lambda y)=Z_{j}-\textsc{O}_{j}$$
where $Z_{j}:=[Z(\lambda j\Delta)-Z(\lambda(j-1)\Delta)]\overset
{d}=Z(\lambda\Delta)$. Furthermore $\textsc{O}_{l}$ appears in
each $\tau_{i}$ for $i\ge l.$ Hence the suggestion by BNS to try
to sample the iid pairs in~\mref{pair2}.

Indeed the joint distribution of the $(\tau_{1},\ldots,\tau_{n})$
is in general complex. However one can easily obtain its joint
Laplace transform. It is with this fact that we argue that the
primary stumbling block which currently prevents one from
integrating out the infinite-dimensional components in the
likelihood, is inherent from the Normal distribution of
$X_{i}|\vartheta,\tau$. Quite simply the Normal assumption yields
exponential terms of the form
$$
{\mbox e}^{-A^{2}_{i}/(2\tau_{i})}{\mbox { rather than }}{\mbox
e}^{-\tau_{i}A^{2}_{i}}.
$$
We will show in the forthcoming sections how to apply a Bessel
integral representation, which does not depend on the distribution
of $(\tau_{1},\ldots,\tau_{n})$ to obtain expressions for
likelihood based on quite general candidates for $\tau.$ First
however we will describe the very remarkable and unique properties
of the OU-$\Gamma$ model in section 3 which does not require this
approach.
\subsection{Some points about GGC and Gamma processes}
We will be making extensive use of the basic elements of the
theory of Generalized Gamma Convolutions(GGC) which can be found
in Bondesson~(1979, 1992) and Thorin~(1977).  GGC are a sub-class
of infinitely divisible random variables. A nice point is that
they all have the important {\it self-decomposability} property.
This has an interesting consequence since it is well known that
$v(t)$ is a stationary OU process if and only if
$v(0)\overset{d}=v(t)$ is self-decomposable. See for instance
Wolfe~(1992), Jurek and Vervaat (1983), Sato~(1999), Jeanblanc,
Pitman and Yor~(2002) or BNS~(2001a, Theorem 1) for a more precise
statement. That is, there is a large subclass of OU models which
all have GGC laws. Some important examples of GGC random variables
and corresponding processes are GIG laws, Stable laws of index
$0<\alpha<1$, and of course Gamma random variables.

Important, from our point of view, is that a random variable is a
GGC if and only if its L\'evy exponent is expressible as \Eq
\int_{0}^{\infty}d_{\theta}(\omega x)\nu(dx)\label{LGGC}\EndEq for
some arbitrary sigma-finite measure satisfying appropriate
conditions so that~\mref{LGGC} is finite and where \Eq
d_{\theta}(\omega)=\theta\log(1+\omega)=\int_{0}^{\infty}(1-{\mbox
e}^{-\omega s})\theta s^{-1}{\mbox e}^{-s}ds.\label{Gamlevy}\EndEq
corresponding to the L\'evy exponent of a Gamma random variable
with shape parameter $\theta.$ That is to say $d_{\theta}$ is a
special case of $\psi$ and moreover the L\'evy density of a
corresponding Gamma process is given by
$$
\rho_{\theta}(ds)=\theta s^{-1}{\mbox e}^{-s}ds{\mbox { for }}s>0
$$
It then follows that the L\'evy density of a GGC is given by \Eq
\theta s^{-1}\int_{0}^{\infty}{\mbox
e}^{-s/r}\nu(dr)\label{LevGGC}.\EndEq As a consequence, if we
denote a Gamma process on a Polish space $\Xcr$ with sigma finite
shape measure $\theta \nu^{*}$ as $G_{\theta \nu^{*}}$, then
\mref{LGGC} is significant as it coincides with the L\'evy
exponent of an arbitrary Gamma process mean functional say
$G_{\theta \nu^{*}}(g)=\int_{\Xcr}g(x)G_{\theta \nu^{*}}(dx)$
where by a change of variable, $R=g(X)$, one can write
equivalently in distribution as $\int_{0}^{\infty}rG_{\theta
\nu}(dr).$

Now we point to a key fact that has not been exploited much in the
literature. First throughout this paper let $T_{\theta}$ denote a
Gamma random variable with shape $\theta$ and scale $1$. Denote
the density of a Gamma random variable with shape $\theta$ and
scale $b>0$ as
$$
\Gcr_{\theta}(y|b)=\frac{b^{-\theta}}{\Gamma(\theta)}y^{\theta-1}{\mbox
e}^{-y/b}{\mbox { for }}y>0.
$$
When $b=1$, we simply write $\Gcr_{\theta}(y).$ Let $(J_{i})$
denote the jump points of a Gamma process and let $(Z_{i})$ denote
the points of a Poisson random measure whose laws are determined
by $\nu,$ which are independent of $(J_{i})$. It is well known
that one can write $G_{\theta
\nu}(dx)=\sum_{i=1}^{n}J_{i}\delta_{Z_i}(dx).$ Furthermore, it
follows that if $Y$ is a GGC random variable then one can always
write
$$
Y\overset{d}=G_{\theta \nu^{*}}(g)\overset {d}=T_{\theta}M_{\theta
\nu}
$$
where $M_{\theta \nu}=\sum_{i=1}^{\infty}(J_{i}/T_{\theta})Z_{i},$
is a random variable independent of $T_{\theta}.$ The independence
property is due to the known fact that the sequence
$(J_{i}/T_{\theta})$ of probabilities is independent of
$T_{\theta}$ which may be written as
$T_{\theta}=\sum_{j=1}^{\infty}J_{i}.$ This property uniquely
characterizes a Gamma process and has nothing do with whether or
not $\nu$ is finite or more generally sigma finite. The sequence
$(J_{i}/T_{\theta})$ is known to have the Poisson-Dirichlet law.

Hence when  $\nu:=H$ is a finite measure, which we will take
without loss of generality to be a probability measure,  a
Dirichlet Process with shape $\theta H$, having total mass $\theta
H(\Xcr)=\theta,$ is defined by the representations $$P_{\theta
H}(dx):=\frac{G_{\theta H}(dx)}{G_{\theta
H}(\Xcr)}=\sum_{j=1}^{\infty}\frac{J_{i}}{T_{\theta}}\delta_{Z_{i}}(dx)$$
where importantly $T_{\theta}=G_{\theta H}(\Xcr)$ is independent
of $P_{\theta H}.$ Setting $\Xcr=[0,\infty)$ one has that $$
M_{\theta H}=\int_{0}^{\infty}xP_{\theta H}(dx)$$ is a Dirichlet
Process mean functional which again is independent of
$T_{\theta}.$ This independent property naturally comes from the
finite dimensional Beta-Gamma calculus based on the classic result
of Lukacs~(1955), which we shall also use. That is, if
$T_{\theta_{i}}$ for $i=1,\ldots, n$ are independent Gamma random
variables with shape $\theta_{i}$ then the sum,
$\sum_{i=1}^{n}T_{\theta_{i}}=T_{\theta*}$, where
$\theta^{*}=\sum_{j=1}^{n}\theta_{i}$ and moreover is independent
of the vector of probabilities $(T_{\theta_{i}}/T_{\theta^{*}})$
which has the Dirichlet distribution with density
$$
D(p_{1},\ldots, p_{n})\propto \prod_{i=1}^{n}p_{i}^{\theta_{i}-1}
$$
where the $\sum_{i=1}^{n}p_{i}=1$. We will denote the fact that a
random vector has Dirichlet law of this type by writing
$\textsc{Dirichlet}_{n}(\theta_{1},\ldots,\theta_{n}).$ Similarly
denote a two parameter beta law as
$\textsc{Beta}(\theta_1,\theta_{2}).$
\subsubsection{Connection to Cifarelli and Regazzini distribution
theory} Because of these observations we are able to exploit the
works of Cifarelli and Regazzini~(1990) and those of subsequent
authors to obtain expressions for the marginal densities of
relevant components of large class of models which we call
OU-FGGC. The FGGC are models with L\'evy density defined
by~\mref{LevGGC} with $\nu:=H.$ This is relevant to both option
pricing and likelihood estimation. We should add that many of
these properties will extend to more general moving average models
where $Z$ is an FGGC BDLP.

The study of properties of Dirichlet process mean functional has
been a major area of interest in Bayesian Nonparametrics. This
line of work was initiated by the paper of Cifarelli and
Regazzini~(1990). One of their important contributions was to
obtain explicit expressions for densities of mean functionals
$M_{\theta H}.$ Let $f_{M_{\theta H}}$ denote the density of
$M_{\theta H}$. Set $H(x)=\int_{0}^{x}H(du).$ Then from Cifarelli
and Regazzini~(1990) or Cifarelli and Melilli~(2000) one has for
$\theta=1$ \Eq f_{M_{H}}(x)=\frac{1}{\pi}\sin(\pi H(x)){\mbox
e}^{-\int_{0}^{\infty}\log(|t-x|)H(dt)} \label{Mdensity}\EndEq and
when $\theta>1$, \Eq f_{M_{\theta
H}}(x)=\frac{\theta-1}{\pi}\int_{0}^{x}{(x-u)}^{\theta-2}\frac{1}{\pi}\sin(\pi
\theta H(u)){\mbox e}^{-\theta\int_{0}^{\infty}\log(|t-u|)H(dt)}du
\label{Mdensity2}\EndEq One can also obtain an expression for the
cdf of $M_{\theta H}$ that holds for all $\theta>0$, we do not
list that here.

\subsubsection{Perfect simulation of $M_{\theta H}$} It is evident
that given the form of the density in $\mref{Mdensity}$ one can in
principle use some sort of rejection sampling procedure to obtain
realizations of $M_{H}.$ With a bit more care one can devise an
efficient method to sample $M_{\theta H}$ for $\theta >1$ using
the density in~\mref{Mdensity2}. Importantly, as pointed out by
Hjort and Ongaro~(2005), when $\theta=m$, where $m=2,3,\ldots$, is
an integer one can use~\mref{Mdensity} to sample $M_{mH}$ based on
the following fact, $$
M_{mH}\overset{d}=\sum_{i=1}^{m}P_{i}M_{1,i}$$ where $(M_{1,i})$
are iid with common distribution equivalent to $M_{H}$ given
by~\mref{Mdensity}. Moreover $(P_{1},\ldots,P_{n})$ is independent
of $(M_{1,i})$ and is $n$-dimensional
$\textsc{Dirichlet}_{n}(1,\ldots,1).$ This can be seen as a simple
consequence of the infinite divisibility of $T_{m}M_{mH}$, where ,
as a consequence, $T_{m}M_{m,H}\overset
{d}=\sum_{i=1}^{m}T_{1,i}M_{1,i}$, and applying the Beta Gamma
calculus. That is further writing $P_{i}=T_{1,i}/T_{m}$, where
$T_{m}=\sum_{j=1}^{m}T_{1,j}$ is independent of $(P_{i}).$ What is
important is that these methods do not rely on the more
computationally burdensome, and otherwise approximate, series
methods. There is however yet another approach which will allow
one to easily \emph{perfectly sample} $M_{\theta H}$ for all
$\theta>0.$

Recently, in the case where $M_{\theta H}$ is almost sure bounded,
Guglielmi, Holmes and Walker~(2002) devise a very simple and
efficient method to obtain perfect samples from the distribution
of $M_{\theta H}$ that works for all $\theta>0.$ We recount the
basic elements of that algorithm. First note that $0\leq a\leq
M_{\theta H}\leq b$ if and only if the support of $H$ is $[a,b]$.
As explained in Guglielmi, Holmes and Walker~(2002), following the
procedure of Propp and Wilson~(1996), one can design an upper and
lower dominating chain starting at some time $-N$ in the past up
to time $0$. The upper chain, say $uM_{\theta H}$, is started at
$uM_{\theta H,-N}=b$, and the lower chain, $lM_{\theta H}$, is
started at $lM_{\theta H,-N}=a$. One runs the Markov chains for
each $n$ based on the equations, \Eq uM_{\theta
H,n+1}=B_{n,\theta}X_{n}+(1-B_{n,\theta}) uM_{\theta H,n}
\label{upper} \EndEq and

\Eq lM_{\theta H,n+1}=B_{n,\theta}X_{n}+(1-B_{n,\theta})
lM_{\theta H,n} \label{lower} \EndEq where the chains are coupled
using the same random independent pairs $(B_{n,\theta}, X_{n})$
where for each $n$, $B_{n,\theta}$ has a Beta$(1,\theta)$
distribution and $X_{n}$ has distribution $H.$ The chains are said
to coalesce when $D=|uM_{\theta H,n}-lM_{\theta H,n}|<\epsilon$
for some small $\epsilon.$ Notice importantly that this method
only requires knowledge of the distribution $H.$

\Remark Vershik, Tsilevich and Yor~(2004) and James~(2005a) are
two examples of applications that directly exploit the
independence property exhibited at the level of the
Gamma/Dirichlet process. See also Diaconis and Kemperman~(1996)
and Diaconis and Freedman~(1999) for more interesting facts.
\EndRemark

\Remark More discussion on the merits of self-decomposability as
it relates to financial applications can be found in Carr, Geman
Madan and Yor~(2005).\EndRemark

\section{Laws and Likelihoods for the
OU-$\Gamma$ model}For $\theta >0$, define a OU-$\Gamma$ process by
setting $Z=G_{\theta}$, where $G_{\theta}$ denote a homogeneous
Gamma process on $[0,\infty)$, i.e. $\nu(dx)=dx$ for
$x\in[0,\infty)$ with law specified by its L\'evy exponent
$d_{\theta}(\omega)$ given in~\mref{Gamlevy}. Letting
$v_{\theta}(t)$ denote the stationary OU-$\Gamma$ it follows that
its L\'evy exponent is \Eq\int_{0}^{1}d_{\theta}(\omega
u)u^{-1}du=\theta \int_{0}^{\infty}(1-{\mbox e}^{-\omega
y})y^{-1}E_{1}(y)dy=-\theta Li_{2}(-\omega)\label{Sv}\EndEq where
$ E_{1}(y)=\int_{y}^{\infty}{\mbox
e}^{-u}u^{-1}du=\int_{1}^{\infty}{\mbox e}^{-uy}u^{-1}du $ is
\emph{Euler's exponential integral}. That is to say the L\'evy
density of $v_{\theta}(0)$ is $\rho_{v_{\theta}}(dy)=\theta
y^{-1}E_{1}(y)dy.$ \Remark In addition to obtaining the form of
the L\'evy density, BNS (2003, p.283) note that the L\'evy
exponent of a OU-$\Gamma$ can be expressed as,
$$
\theta\sum_{j=1}^{\infty}(-1)^{j}\frac{\omega^{j}}{j^{2}}{\mbox {
for }}0\leq \omega<1$$ but they don't equate this with the
dilogarithm function. \EndRemark

The previous discussion indicates that one can implement both
option pricing and likelihood analysis if one can sample the
special case of ~\mref{pair2} given by
$$(G_{\theta}(\lambda \Delta),{\mbox e}^{-\lambda
\Delta}\int_{0}^{\Delta}{\mbox e}^{\lambda y}G_{\theta}(d\lambda
y)).$$ The L\'evy exponent of the second term is given by
$$
\int_{{\mbox e}^{-a}}^{1}d_{\theta}(\omega u)u^{-1}du=\int_{{\mbox
e}^{-a}}^{1}d_{\theta a}(\omega u)F_{a}(du)=-\theta
[Li_{2}(-\omega)-Li_{2}(-\omega e^{-a}))]$$ where \Eq
F_{a}(y)=\int_{{\mbox
e}^{-a}}^{y}\frac{1}{au}du=\frac{\log(y)+\a}{a} \label{cdf}\EndEq
is a cdf for ${\mbox e}^{-a}\leq y\leq 1$. However due to the fact
$G_{\theta}(\lambda \Delta)\overset {d}=G_{\theta a F_{a}}([{\mbox
e}^{-a},1])\overset{d}=T_{\theta a}$ this is equivalent to
sampling the pair $$(T_{\theta a},\int_{{\mbox
e}^{-a}}^{1}xP_{\theta a F_{a}}(dx))$$ where for ${\mbox
e}^{-a}\leq y\leq 1$
$$
P_{\theta a F_{a}}(dy)=\frac{G_{\theta a F_{a}}(d y)}{T_{\theta
a}}
$$
is a Dirichlet process random probability measure with shape
parameter $\theta a F_{a}.$

\Remark We shall use the notation $M_{\theta a}$ rather than the
perhaps more accurate $M_{\theta a F_{a}}$ where it is understood
that $F_{a}$ is defined in~\mref{cdf}\EndRemark

We discuss some of the implications of these facts in the next two
propositions.
\begin{prop}For each fixed $\Delta>0$, and $a=\lambda \Delta$ set
$Y_{\theta a}:={\mbox e}^{-a}\int_{0}^{\Delta}{\mbox e}^{\lambda
y}G_{\theta}(d\lambda y)$, where $G_{\theta}$ is a homogeneous
Gamma process. It follows that
$G_{\theta}(a)=\int_{0}^{t}G_{\theta}(d\lambda
y)\overset{d}=T_{\theta a}$. Additionally, the following
distributional properties hold. \Enumerate
\item[(i)]Let $M_{\theta a}:=\int_{0}^{1}xP_{\theta a F_{a}}(dx)$
denote a Dirichlet process mean functional based on the shape
parameter $\theta a F_{a}.$ Then for each fixed $\Delta$, one has
the coordinate-wise equivalence in joint distribution,
$$
(G_{\theta }(a), Y_{\theta a})\overset{d}=(T_{\theta a}, T_{\theta
a} M_{\theta a})
$$
where $M_{\theta a}$ is independent of $T_{\theta a}.$ Furthermore
${\mbox e}^{-a}\leq M_{\theta a}\leq 1$ almost surely.
\item[(ii)]$\int_{0}^{t}(1-{\mbox e}^{-\lambda
(\Delta-y))}G_{\theta}(d\lambda y)=G_{\theta}(a)-Y_{\theta
a}\overset {d}=T_{\theta a}[1-M_{\theta a}]$
\item[(iii)] $(G_{\theta}(a)-Y_{\theta a},Y_{\theta a})\overset{d}=(T_{\theta a}[1-M_{\theta a}],T_{\theta a}M_{\theta
a})$ \EndEnumerate\qed
\end{prop}
\Proof The result is already established by our construction and
appealing to the unique independence property of the
Gamma/Dirichlet process. However since the joint equivalence in
statement (i) is the key factor separating the OU-$\Gamma$ from
other OU processes, hence quite delicate, we will check it via
joint Laplace transforms. Evaluating the joint Laplace transform
of the $(G_{\theta }(a), Y_{\theta a})$ at points
$(\omega_{1},\omega_{2})$, it is easily seen that the joint L\'evy
exponent is $$\int_{{\mbox e}^{-a}}^{1}d_{\theta
a}(\omega_{1}+\omega_{2}u)F_{a}(du).
$$  Now being careful to use only the independence property of
$T_{\theta a}$ and $M_{\theta a}$ and the fact that $M_{\theta a}$
is a Dirichlet process mean functional we proceed as follows.
Write $\omega_{1}T_{\theta a}+\omega_{2}T_{\theta a} M_{\theta
a}=T_{\theta a}[\omega_{1}+\omega_{2}M_{\theta a}]:=W.$
Furthermore note the $\omega_{1}+\omega_{2}M_{\theta
a}=\int_{0}^{1}(\omega_{1}+\omega_{2}x)P_{\theta a
F_{a}}(dx)=P_{\theta aF_{a}}(g)$, for
$g(x)=\omega_{1}+\omega_{2}x.$ Now by independence of $T_{\theta
a}$ and $M_{\theta a}$ the joint Laplace transform, taking
expectation with respect to the Gamma law first is,
$$
\E[{\mbox e}^{-W}]=\E[{(1+\omega_{1}+\omega_{2}\omega M_{\theta
a})}^{-\theta a}]=\E[{(1+P_{\theta aFa}(g))}^{-\theta a}]
$$
Now appealing to the well-known identity of Cifarelli and
Regazzini~(1990) it follows that
$$
\E[{(1+P_{\theta aFa}(g))}^{-\theta a}]=\E[{\mbox e}^{-G_{\theta a
F_{a}}(g)}]
$$
which is the desired result. The above argument indeed establishes
the proof but the very special nature of the result perhaps will
not be fully clear until one reads section 3.8.\EndProof

The next result describes the distribution of $v_{\theta}(0)$ in
the stationary case.

\begin{prop} Let $v_{\theta}(0)$ have distribution described by the L\'evy exponent~\mref{Sv}. Let
$G_{\theta \nu}$ denote a (non-finite) Gamma process on $[0,1]$
with $\nu(du)=u^{-1}du$ where $\int_{0}^{1}\nu(du)=\infty.$ Then
$v_{\theta}(0)$ is a generalized Gamma convolution (GGC) such that
$$
v_{\theta}(0)\overset{d}=\int_{0}^{1}xG_{\theta
\nu}(dx)\overset{d}=T_{\theta}{\tilde M}_{\theta},
$$
where ${\tilde M}_{\theta}=M_{\theta \nu}$ is independent of
$T_{\theta}$ but is not a Dirichlet process mean functional.
Furthermore, for each fixed $\theta$, the distribution of ${\tilde
M}_{\theta }$ is characterized by its generalized Cauchy-Stieltjes
transform,
$$
\E[{\mbox e}^{-\omega v_{\theta}(0)}]=\E[{(1+\omega {\tilde
M}_{\theta })}^{-\theta}]={\mbox e}^{\theta Li_{2}(-\omega)}
$$\qed
\end{prop}

\Remark It is quite possible to obtain an explicit form of the
density of $v_{\theta}(0)$ by using standard inversion results for
characteristic functions and noting the relationship of the
complex valued dilogarithm function to the \emph{Inverse Tangent
Integral},
$$
Ti_{2}(y)=\int_{0}^{y}\frac{\arctan(u)}{u}du,
$$
which is the imaginary part of the complex valued dilogarithm
function, and \emph{Clausen's Function}. For more details see
Maximon~(2003). \EndRemark

Recapping, Proposition 3.1 shows that the distribution of
$(G_{\theta}(a), Y_{\theta a})$ is determined by the distribution
of the independent random variables $(T_{\theta a}, M_{\theta
a}).$  Among OU processes discussed here, the independence
property is unique to OU-$\Gamma$ processes. Additionally, as we
shall see this pair may be sampled exactly due to the fact that
$M_{\theta a F_{a}}$ is a Dirichlet process mean functional. On
the other hand Proposition 3.2 shows that although
$v_{\theta}(0)\overset{d}=T_{\theta}{\tilde M}_{\theta }$ is a
GGC, the results for the Dirichlet process do not apply to
${\tilde M}_{\theta }$ and we otherwise do not have a tractable
expression for the explicit density of $v_{\theta}(0)$. However,
we do believe that a careful use of the relationships mentioned in
Remark 6 will lead to an explicit form. The next proposition,
using the work of Cifarelli and Regazzini~(1990), provides more
details for the distribution of $M_{\theta a }$ and shows also
that one can use the Dirichlet process results to obtain a good
approximate for the distribution of $v_{\theta}(0).$

\begin{prop} For each $0<a=\lambda \Delta<\infty$ and $\theta>0$, let
$Y_{\theta a}\overset{d}={\mbox e}^{-\lambda
\Delta}\int_{0}^{\Delta}{\mbox e}^{\lambda y}G_{\theta}(d\lambda
y)$ denote an infinitely divisible random variable with L\'evy
exponent,
$$
\int_{{\mbox e}^{-a}}^{1}d_{\theta}(\omega u)u^{-1}du=\int_{{\mbox
e}^{-a}}^{1}d_{\theta a}(\omega u)F_{a}(du)=-\theta
[Li_{2}(-\omega)-Li_{2}(-\omega e^{-a}))]
$$
where $F_{a}(du)=a^{-1}u^{-1}du$ is the density of a random
variable taking its values in the interval $[{\mbox e}^{-a},1].$
Then the following results hold \Enumerate
\item[(i)] $Y_{\theta a}\overset{d}=T_{\theta a}M_{\theta
a}$, where $M_{\theta a}=\int_{{\mbox e}^{-a}}^{1}xP_{\theta
aF_{a}}(dx)$ is a Dirichlet process mean functional.
\item[(ii)]The L\'evy density of $Y_{\theta a}$ is
$\rho_{\theta a}(dy)=\theta y^{-1}[E_{1}(y)-E_{1}(ye^{a})]dy.$
Hence the cumulants of $Y_{\theta a}$ are for each integer $j$,
$$
\theta\int_{0}^{\infty}y^{j-1}[E_{1}(y)-E_{1}(ye^{a})]dy=\theta
\frac{\Gamma(j)}{j}(1-{\mbox e}^{-aj})
$$
\item[(iii)] When $\theta a=1$, the
density of $M_{1}$ is given by\Eq
\frac{1}{\pi}\sin\(\[\frac{-\pi\log(x)}{a}\]\)
x^{\frac{1}{a}[1-\log(x)]-1}{\mbox e}^{\frac{\pi^{2}}{3a}}{\mbox
e}^{\frac{-1}{a}[Li_{2}(x)+Li_{2}(\frac{{\mbox e}^{-a}}{x})]},
\label{dilogden}\EndEq for ${\mbox e}^{-a}\leq x\leq 1.$
\item[(iv)]When $\theta a=1$, the
density of $V_{a}:=-\log(M_{1})/a$ is given by\Eq
\frac{1}{\pi}\sin(\pi v) {\mbox e}^{-[v+v^{2}]}{\mbox
e}^{\frac{\pi^{2}}{3a}}{\mbox e}^{\frac{-1}{a}[Li_{2}({\mbox
e}^{-av})+Li_{2}({\mbox e}^{-a(1-v)})]}, \label{logdilogden}\EndEq
for $0\leq v\leq 1.$
\item[(v)]When $\theta a>1$, the
density of $M_{\theta a}/a$ is given by\Eq \frac{\theta
a-1}{\pi}\int_{-\log x/a}^{1}{(x-{\mbox e}^{-va})}^{\theta
a-2}\sin(\pi \theta av) {\mbox e}^{-\theta a[v+v^{2}]}{\mbox
e}^{\frac{\theta \pi^{2}}{3}}{\mbox e}^{-\theta [Li_{2}({\mbox
e}^{-av})+Li_{2}({\mbox e}^{-a(1-v)})]},
\label{2logdilogden}\EndEq for $0\leq v\leq 1.$
\item[(vi)]
If $\theta a=m$, where $m=2,3,\ldots$ is an integer, then
$M_{m}\overset{d}=\sum_{i=1}^{N}W_{i}M_{1,i}$, where $(M_{1,i})$
are iid with density~\mref{dilogden} and independent of
$(M_{1,i})$, $(W_{i}=T_{1,i}/\sum_{j=1}^{N}T_{1,j})$, where ,
$T_{1,i}\overset{d}=T_{1}$ are iid, is a $\textsc{Dirichlet}_{n}
(1,\ldots,1)$ $n$-dimensional vector.

\item[(vii)]$Y_{\theta a}$ converges in distribution to
$v_{\theta}(0)$ as ${\mbox e}^{-a}\rightarrow 0.$\EndEnumerate\qed
\end{prop}
\Proof Most of the results are immediate from our previous
discussioms. The forms of the density arises from application of
Cifarelli and Regazzini~(1990) which amounts to explicitly
calculating $\int_{{\mbox e}^{-a}}^{1}\log(|t-x|)F_{a}(dt)$
expressed in terms of the dilogarithm function.\EndProof The last
result in this section gives a completely tractable description of
the conditional distribution of the log asset price at time $t$
given information up to time $s$. This type of result is pertinent
to option pricing as discussed in BNS(2001a, 6.2) and Nicolato and
Vernados~(2003).

\begin{prop} Let $x^{*}_{\theta}(t)$ be defined as in~\mref{BNS} with
$Z=G_{\theta}$. Additionally for $0\leq s<t$, set $\Delta=(t-s)$
and define $h(\Delta,s)=(1-{\mbox e}^{-\lambda
\Delta})v_{\theta}(s)$ and
$\mu^{*}_{s}=\mu\Delta+x^{*}_{\theta}(s)+\beta h(\Delta,s).$ Then
the conditional density of
$x^{*}_{\theta}(t)|x^{*}_{\theta}(s),v_{\theta}(s)$ is given by
$$
\int_{0}^{\infty}\phi(x|\mu^{*}_{s}+\beta y,h(\Delta,s)+
y)q_{\theta a}(y)dy
$$
where $q_{\theta a}(y)=\int_{{\mbox e}^{-a}}^{1}\Gcr_{\theta
a}(y|(1-v))f_{M_{\theta a}}(v)dv.$ When $\theta a=1,$ $$
q_{1}(y)=\int_{0}^{1}\Gcr_{1}(y|(1-{\mbox e}^{-v}))
f_{V_{a}}(v)dv$$ where $f_{V_{a}}$ is the density of $V_{a}$ given
in~\mref{logdilogden}.\qed
\end{prop}

\subsection{Perfect simulation of $M_{\theta a}$}
Due to the fact that the dilogarithm function, $Li_{2}(x)$, is a
well-understood special function, which is available in many
computational packages, it is evident that the densities in
$\mref{dilogden}$ and $\mref{logdilogden}$ can be exactly sampled
using a rejection procedure. Again based on the discussion in
section 2.3.2 Statement (vi) of Proposition 3.3 shows that one can
use this fact to easily obtain samples of $M_{m}$, and hence
$Y_{m}$, for any integer $m.$ With a bit more care one can devise
an efficient method to sample $M_{\theta a}$ for $\theta a>1$
using the density in~\mref{2logdilogden}. One can also use the
perfect sampling method described in 2.3.2 for all $\theta a$,
based on $uM_{\theta a,-N}=1$ and $lM_{\theta a,-N}={\mbox
e}^{-a}$, $B_{n,\theta a}$ is \textsc{Beta} $(1,\theta a)$ and
$X_{n}$ has distribution $F_{a}$

\subsection{BNS OU-$\Gamma$ likelihood inference}
The results in the previous section now give the ingredients to
perform likelihood based statistical inference via  simple exact
sampling. Here we describe a bit more about the distribution of
$\tau_{i}$ in the OU-$\Gamma$ case and then extend the discussion
to randomly sampled times.
\begin{prop}Define for $\Delta>0$, $a=\lambda \Delta$ and $i=1,\ldots, n$
$\tau_{\theta,i}:=\tau_{\theta}(i\Delta)-\tau_{\theta}((i-1)\Delta)$,
by setting $Z=G_{\theta}$ in  ~\mref{tdiff}. Furthermore, let
$r_{i}={\mbox e}^{\lambda i\Delta}$ for $i=1,\ldots, n$, with
$r_{0}=1$. Then it follows that, for $i=1,\ldots, n$, \Eq
\lambda\tau_{\theta,i}=(1-{\mbox e}^{-\lambda
\Delta})v_{\theta}((i-1)\Delta)+T_{i}[1-M_{i}] \label{tausimp}
\EndEq with,
$$
v_{\theta}((i-1)\Delta)={\mbox e}^{-\lambda
(i-1)\Delta}[v_{\theta}(0)+\sum_{j=1}^{i-1}r_{j}T_{j}M_{j}]
$$
where $(T_{i},M_{i})$ are iid pairs independent of
$v_{\theta}(0).$ Additionally, for each fixed $i$, $T_{i}$ and
$M_{i}$ are independent with distributions specified by
$T_{i}\overset{d}=T_{\theta a}$ and $M_{i}\overset{d}=M_{\theta
a}.$ This implies that likelihood inference for the
model~\mref{BNSlik} may be obtained from the joint distribution of
$(X_{i},T_{i},M_{i},v_{\theta}(0))$ given by \Eq
\[\prod_{i=1}^{n}\phi(X_{i}|\mu\Delta+\beta
\tau_{\theta,i},\tau_{\theta,i})\Gcr_{\theta
a}(t_{i})f_{M}(v_{i})\]f_{v_{\theta}(0)}(w) \label{auglik}\EndEq
where $\tau_{\theta,i}$ is expressed as in~\mref{tausimp}, with
$T_{i}=t_{i}$,$M_{i}=v_{i}$, and $v_{\theta}(0)=w$. \qed
\end{prop}
A Bayesian procedure, which involves placing a prior on
$\vartheta=(\mu,\beta,\lambda,\varrho)$, is quite natural and
otherwise proceeds by standard arguments, in this setting. That is
letting $\pi(\vartheta)$ denote a prior joint density it follows
that a posterior distribution of $\vartheta|\X$ is determined by a
posterior distribution of
$\vartheta,(T_{i},M_{i}),v_{\theta}(0)|\X$ which is proportional
to
$$
\pi(\vartheta)\[\prod_{i=1}^{n}\phi(X_{i}|\mu\Delta+\beta
\tau_{\theta,i},\tau_{\theta,i})\Gcr_{\theta
a}(t_{i})f_{M}(v_{i})\]f_{v_{\theta}(0)}(w)
$$
\Remark The likelihood in~\mref{BNSlik} for the OU-$\Gamma$ case
obviously is obtained by integrating out the pertinent independent
quantities in~\mref{auglik}. Due to the Gamma distributions, the
answer could be expressed in terms of integrals with respect to
modified Bessel functions. Or otherwise a subclass of Generalized
Inverse Gaussian(GIG) random variables.\EndRemark \Remark Note
that in practice we can approximate a draw from the distribution
of $v_{\theta}(0)$ by using instead $Y_{\theta \delta}$ for
${\mbox e}^{-\delta}$ small. Otherwise, if strict stationarity
$v_{\theta}(t)$is not a concern, one can certainly use any
positive distribution for $v_{\theta}(0).$ \EndRemark

\subsection{The likelihood via a connection to Variance Gamma
processes} Recall in the stationary case that according to
Proposition 3.2. $v_{\theta}(0)\overset{d}=T_{\theta}{\tilde
M}_{\theta}$, where ${\tilde M}_{\theta }$ is not a Dirichlet
process mean functional. However this point allows one to write
$\tau_{\theta}$ and $(\tau_{\theta,1},\ldots, \tau_{\theta,n})$ in
terms of a product of a Gamma random variable and another
independent random variable. Specifically, for $a=\lambda \Delta$,
one may write
$$
\tau_{\theta,i}=T_{\theta(1+na)}S_{i}
$$
where for $i=1,\ldots, n$
$$
\lambda S_{i}=(1-{\mbox e}^{-a}){\mbox e}^{-a
(i-1)}\[\frac{T_{\theta}}{T_{\theta
(1+na)}}\tilde{M}_{\theta}+\sum_{j=1}^{i-1}\frac{T_{j}}{T_{\theta(1+na)}}r_{j}M_{j}\]+\frac{T_{i}}{T_{\theta(1+na)}}
[1-M_{i}].
$$
The vector $\S=(S_{1},\ldots, S_{n})$ is independent of
$T_{\theta(1+na)}$ which can be written as
$T_{\theta}+\sum_{i=1}^{n}T_{i}$. We may also write
$$
\lambda S_{i}=(1-{\mbox e}^{-a}){\mbox e}^{-a (i-1)}\[P_{
n+1}\tilde{M}_{\theta}+\sum_{j=1}^{i-1}P_{j}r_{j}M_{j}\]+P_{i}[1-M_{i}].
$$
where $P_{n+1}=1-\sum_{j=1}^{n}P_{j},$  and $(P_{1},\ldots,
P_{n+1})$ is $\textsc{Dirichlet}_{n+1}(\theta a,\ldots,\theta
a,\theta)$ independent of all other random variables. Recall now
that a $GIG(\nu,\delta, \gamma)$ random variable has density given
by
$$
g(x|\nu,\delta,
\gamma)=\frac{{(\gamma/\delta)}^{\nu}}{2K_{\nu}(\delta\gamma)}x^{\nu-1}{\mbox
e}^{-\frac{1}{2}(\delta^{2}x^{-1}+\gamma^{2}x)}{\mbox { for }}x>0
$$
where $K_{\nu}$ is a modified Bessel function. Recall also that
$K_{\nu}(x)=K_{-\nu}(x)$.

Additionally we will exploit the following nice feature of
$K_{\nu}(x).$ Suppose that for a $m=0,1,2\ldots$, the
$|\nu|=m+1/2$, where $|\nu|$ denotes absolute value, then we can
use the fact that \Eq K_{m+1/2}(x)=\sqrt{\frac{\pi}{2x}}{\mbox
e}^{-x}\sum_{k=0}^{m}\frac{(m+k)!}{k!(m-k)!2^{k}}x^{-k}
\label{bessimp} \EndEq See for instance Pitman~(1999, eq. (40))
for a probabilistic interpretation of~\mref{bessimp}.

This facts leads to the following description of the likelihood.
\begin{thm} The observations according to~\mref{BNSx} can be
represented as $X_{i}=\mu\Delta+\beta
T_{\theta(1+na)}S_{i}+{[T_{\theta(1+na)}S_{i}]}^{1/2}\epsilon_{i},
$ in the OU-$\Gamma$ case. Setting
$\gamma^{2}=[2+\beta^{2}\sum_{j=1}^{n}S_{i}]$ and
$\delta^{2}=\sum_{j=1}^{n}A^{2}_{j}/(2S_{j})$,
$\kappa=\theta(1+na)$ and $\nu=\kappa-n/2,$ and $a=\lambda\Delta.$
The following results hold. \Enumerate
\item[(i)]The likelihood in~\mref{BNSlik} can be written as,
$$
\Lcr(\X|\vartheta)={\mbox
e}^{n\bar{A}\beta}\E_{\vartheta}\[\frac{2K_{\nu}(\delta\gamma)}{{(\gamma/\delta)}^{\nu}\Gamma(\kappa)}\prod_{i=1}^{n}\frac{1}
{\sqrt{2\pi}}S^{-1/2}_{i}\]
$$
\item[(ii)]If $\theta$ and $a$ are chosen such $|\nu|=m+1/2$, for
$m=0,1,2,\ldots$, then
$$
\Lcr(\X|\vartheta)={\mbox
e}^{n\bar{A}\beta}\sum_{k=0}^{m}\frac{(m+k)!}{k!(m-k)!2^{k}}
\E_{\vartheta}\[{\mbox
e}^{-\delta\gamma}\frac{2{(\delta\gamma)}^{-k}\gamma^{-1}}{{(\gamma/\delta)}^{m}\Gamma(\kappa)}\prod_{i=1}^{n}\frac{1}
{\sqrt{2\pi}}S^{-1/2}_{i}\]\sqrt{\frac{\pi}{2}}
$$
As a special case $|\nu|=m+1/2$ for all $n$, if $\theta a=1/2$ and
$\theta=m+1/2.$
\item[(iii)]If additionally $m=0$, that is $\theta=1/2$ and $a=1$, then
$$
\Lcr(\X|\vartheta)={\mbox e}^{n\bar{A}\beta}
\E_{\vartheta}\[{\mbox
e}^{-\delta\gamma}\frac{2\gamma^{-1}}{\Gamma(\kappa)}\prod_{i=1}^{n}\frac{1}
{\sqrt{2\pi}}S^{-1/2}_{i}\]\sqrt{\frac{\pi}{2}}
$$
\EndEnumerate In all cases the distribution of $(S_{1},\ldots,
S_{n})$ is completely determined by the $2n+2$ independent random
variables with joint density $\[\prod_{i=}^{n}\Gcr_{\theta
a}(t_{i})f_{M_{\theta a}}(v_{i})\]\Gcr_{\theta}(t)f_{\tilde
M_{\theta}}(w)$\qed
\end{thm}
The next result in effect serves to make clear Theorem 3.1 but
also highlights the possibility, from a practical point of view,
for more data augmentation procedures
\begin{prop} Consider the setup and notation in Theorem 3.1. Additionally
define $\beta^{2}_{*}=\beta^{2}\sum_{j=1}^{n}S_{i}.$ Then it is
clear that
$$
\frac{2K_{\nu}(\delta\gamma)}{{(\gamma/\delta)}^{\nu}\Gamma(\kappa)}=
\int_{0}^{\infty}y^{-\frac{n}{2}}{\mbox
e}^{-\frac{1}{2}(\delta^{2}y^{-1}+\beta^{2}_{*}y)}\Gcr_{\kappa}(y)dy,
$$
which leads to another expression of the likelihood
$\Lcr(\X|\vartheta)$. Thus statistical inference may be based on
simulation from the joint density
$$
\Gcr_{\kappa}(y)\[\prod_{i=}^{n}\Gcr_{\theta a}(t_{i})f_{M_{\theta
a}}(v_{i})\]\Gcr_{\theta}(t)f_{\tilde M_{\theta}}(w).$$ Based on
this fact one has that if a random variable $V$ has the density
$\Gcr{\kappa}$ relative to the representation of
$\Lcr(\X|\vartheta)$ then the posterior distribution of
$V|\S,\X,\vartheta$ is $GIG(\nu,\delta, \gamma)$ with parameters
specified by Theorem 3.1. and Proposition 3.5
\end{prop}

\Remark One notes that the expressions in statements (ii) and
(iii) of Theorem 3.1 are quite manageable. Here one is perhaps
taking the view that $\theta$ and $\lambda \Delta$ are chosen to
ease computations. However note that in statement (ii) that $m$,
whose parameter space is $\{0,1,\ldots, \}$ becomes a viable and
flexible parameter of interest from a modelling point of view. The
expression in statement (i) is also quite amenable to Monte-Carlo
estimation approaches. \EndRemark \Remark By Variance Gamma
processes we are loosely referring to the work of Madan, Carr and
Chang~(1998), see also Carr, Geman, Madan, and Yor~(2003). It is
evident that all OU-GGC models exhibit similar properties. That is
if the BDLP $Z$ is a GGC then analogues of Theorem 3.1 and
Proposition 3.6 have exactly the same form. However, in contrast
to the OU-$\Gamma$ case, one still does not have an obvious way to
sample from the distribution of $\S.$ \EndRemark

\subsection{Bayesian estimation and related comments} We have shown that the distribution of
$(\tau_{\theta,1},\ldots, \tau_{\theta,n})$ is determined by
$2n+2$ independent random variables whose distributions can be
perfectly sampled or in the case of $v_{\theta}(0)$ approximated
with arbitrary accuracy. We also note that the explicit densities
that we have given for $M_{\theta a}$ definitely have practical
utility, whereby rejection methods can be used. We also believe
they are interesting from a mathematical point of view as they may
have connections to application in physics or analytic
combinatorics. These are places where the dilogarithm function
appears often. However, in terms of practical simplicity it is
perhaps easier to use the perfect simulation schemes which work
for all values of $\theta a$ and only require simulation from beta
random variables and the distribution $F_{a}.$  Also, in regards
to $v_{\theta}(0)$, we note again that in the case of not strictly
stationary OU-$\Gamma$ models, we may choose $v_{\theta}(0)$ to
have any distribution. However Theorem 3.1 suggests there are some
quite interesting simplifications that occur if we choose
$v_{\theta}(0)=T_{\theta}W$, where $W$ denotes a random variable
independent of $T_{\theta}.$ We note again that all GGC random
variables have this form including the class of GIG models.

Armed with the information that we have provided one can construct
a variety of efficient simulation based techniques. Here we
briefly highlight the Bayesian approach. Primarily this is due to
the fact that a Bayesian approach is essentially an approach
involving integration and hence is a quite natural for Monte-Carlo
based estimation. It is in many respects quite similar to
Bootstrap techniques. We now mention some well known points about
Bayesian estimation. Suppose that $\pi(\vartheta)$ is a prior
distribution of the unknown parameters. Then, as is well known,
the fundamental object of interest is to obtain the posterior
distribution of $\theta|\X$, which is given by
$$
\pi(\vartheta|\X)\propto \pi(\vartheta)\Lcr(\X|\vartheta)
$$
Estimation of some parameter $h(\vartheta)$ can then be cast in
terms of integration, \Eq
E[h(\vartheta)|\X]=\int_{\Theta}h(u)\pi(du|\X)=\frac{\E\[h(\vartheta){\mbox
e}^{n\bar{A}\beta}\frac{2K_{\nu}(\delta\gamma)}{{(\gamma/\delta)}^{\nu}\Gamma(\kappa)}\prod_{i=1}^{n}\frac{1}
{\sqrt{2\pi}}S^{-1/2}_{i}\]}{\E\[{\mbox
e}^{n\bar{A}\beta}\frac{2K_{\nu}(\delta\gamma)}{{(\gamma/\delta)}^{\nu}\Gamma(\kappa)}\prod_{i=1}^{n}\frac{1}
{\sqrt{2\pi}}S^{-1/2}_{i}\]}\label{Bayes1}\EndEq where the
denominator should be understood as,
$$
\Lcr(\X)=\int_{\Theta}\E_{\vartheta}\[{\mbox
e}^{n\bar{A}\beta}\frac{2K_{\nu}(\delta\gamma)}{{(\gamma/\delta)}^{\nu}\Gamma(\kappa)}\prod_{i=1}^{n}\frac{1}
{\sqrt{2\pi}}S^{-1/2}_{i}\]\pi(d\vartheta).
$$
For instance, the posterior probability that $\vartheta$ is in
some region $B$ can be evaluated by choosing $h(x)=I\{x\in B\}.$
Since Bessel functions, such as $K_{\nu}(x),$ are available in
standard mathematical computer packages, one can just draw from
the joint distribution of $(\vartheta,\S)$, which is readily
available from our results. That is for $l=1,\ldots, B$ draw iid
random vectors $(\vartheta_{l},S_{1,l},\ldots,S_{n,l})$
then~\mref{Bayes1} is approximated by \Eq
\frac{\sum_{l=1}^{B}h(\vartheta_{l}){\mbox
e}^{n\bar{A}\beta_{l}}\frac{2K_{\nu_{l}}(\delta_{l}\gamma_{l})}{{(\gamma_{l}/\delta_{l})}^{\nu}\Gamma(\kappa_{l})}\prod_{i=1}^{n}\frac{1}
{\sqrt{2\pi}}S^{-1/2}_{i,l}} {\sum_{l=1}^{B}{\mbox
e}^{n\bar{A}\beta_{l}}\frac{2K_{\nu_{l}}(\delta_{l}\gamma_{l})}{{(\gamma_{l}/\delta_{l})}^{\nu}\Gamma(\kappa_{l})}\prod_{i=1}^{n}\frac{1}
{\sqrt{2\pi}}S^{-1/2}_{i,l}} \label{simpest}.\EndEq The nice
feature of basic iid Monte-Carlo type estimator
like~\mref{simpest} is that accuracy issues are well-understood
and are less dependent on the sample size. Here accuracy increases
as $B$ increases.

One can of course develop more sophisticated importance sampling
and MCMC methods based on well-known ideas. These may involve
sampling from the posterior distributions. For instance, our
results show that the posterior distribution of $\vartheta|\X$ can
be obtained by working with the posterior distributions of
$\vartheta|\X,\S,V$ and $\S,V|\X,\vartheta$ where for instance
$V|\S,\X,\vartheta$ has a $GIG(\nu,\delta, \gamma)$ distribution.
All other conditionals can be easily deduced by various
augmentations of the expressions given in Theorem 3.1 and
Proposition 3.6.
\subsection{OU-$\Gamma$ processes with possibly random scale parameter}
Up to this point we have assumed that $G_{\theta}$ was a
homogeneous Gamma process with scale parameter equal to $1$. This
was done mainly for notational convenience. However, it follows
from our analysis that the introduction of a scale parameter say
$\zeta$ can be used as a powerful modeling tool. Naturally a scale
parameter can just be introduced by replacing $G_{\theta}$ with
$\zeta G_{\theta}$ throughout. However an important fact is that
if we use $\zeta G_{\theta}$, the vector $\S$ described in section
3.3 still does not depend on $\zeta$. This means that one can now
write
$$
X_{i}=\mu\Delta+\beta \zeta T_{\theta(1+na)}S_{i}+{[\zeta
T_{\theta(1+na)}S_{i}]}^{1/2}\epsilon_{i}.$$ Note that if $\zeta$
is fixed then all our results carry over without change. This
means extending the model to the case where $\zeta$ is random is
straightforward. The main feature being that we would now be
working with a Gamma scale mixture, based on $\zeta
T_{\theta(1+na)}$, which can be used to introduce more
distributional modeling flexibility.

\subsection{Likelihoods for Superpositioned OU-$\Gamma$}BNS~(2001a, p.178) propose the idea of
superpositions of independent OU processes to alter the
auto-correlation structure. Here, letting $p$ denote a positive
integer, and $(w_{1},\ldots,w_{p})$ a possibly unknown vector of
positive terms summing to $1$,we discuss briefly a generalization
of Theorem 3.1 to the case where one starts with a superposition
process $v(t|p)=\sum_{j=1}^{p}w_{j}v_{\theta_{j}}(t)$ where for
$j=1,\ldots, p$, $v_{\theta_{j}}(t)$ are independent OU-$\Gamma$
processes which are based on parameters
$(\lambda_{j},\theta_{j})$, in place of $(\lambda, \theta)$.
Obviously the distributional results we have developed apply to
each of the independent components. One uses for instance
$a_{j}=\lambda_{j}\Delta$ and $\theta_{j}a_{j}$ in place of $a$
and $\theta a.$

Let $\tau(t|p)=\sum_{j=1}^{p}w_{j}\tau_{\theta_{j}}(t)$, denote
the integrated volatility where each
$\tau_{\theta_{j}}(t)=\int_{0}^{t}v_{\theta_{j}}(s)ds.$
Additionally the analog of $(\tau_{\theta,1},\ldots,
\tau_{\theta,n})$ is $\tau_{i}:=\tau(i\Delta)-\tau_((i-1)\Delta)$.
Then by similar arguments to the previous section one can write
for $\xi_{n}=\sum_{j=1}^{p}\theta_{j}(1+n\lambda_{j}\Delta),$
$$
\tau_{i}=T_{\xi_{n}}S_{i,p}
$$
where $S_{i,p}:=\tau_{i}/T_{\xi_{n}}$ has an obvious description
by applying our previous results to each component
$\tau_{\theta_{j}}$, and the vector $(S_{1,p},\ldots S_{n,p})$ is
independent of $T_{\xi_{n}}.$
\begin{prop} Let $X_{i}=\mu\Delta+\beta
T_{\xi_{n}}S_{i,p}+{[T_{\xi_{n}}S_{i,p}]}^{1/2}\epsilon_{i}.$ with
terms defined in this section. Set
$\gamma^{2}=[2+\beta^{2}\sum_{j=1}^{n}S_{i,p}]$ and
$\delta^{2}=\sum_{j=1}^{n}A^{2}_{j}/(2S_{j,p})$,
$\kappa=\xi_{n}=\sum_{j=1}^{p}\theta_{j}(1+n\lambda_{j}\Delta)$
and $\nu=\kappa-n/2.$ Let $\vartheta_{p}$ denote the enlarged
parameter space containing unknown quantities such as
$(w_{1},\ldots, w_{p})$, then the likelihood
$\Lcr(\X|\vartheta_{p})$ has the same form as the likelihood in
Theorem 3.1 with appropriate substitutions of the above parameters
and $(S_{1,p},\ldots, S_{n,p})$ in place of $(S_{1},\ldots,
S_{n}).$ In particular, \Enumerate
\item[(i)]if $[\sum_{j=1}^{p}\theta_{j}\lambda_{j}]\Delta=1/2$ and $
\sum_{j=1}^{p}\theta_{j}=m+1/2$ for $m=0,1,2,\ldots$, then
$$
\Lcr(\X|\vartheta_{p})={\mbox
e}^{n\bar{A}\beta}\sum_{k=0}^{m}\frac{(m+k)!}{k!(m-k)!2^{k}}
\E\[{\mbox
e}^{-\delta\gamma}\frac{2{(\delta\gamma)}^{-k}\gamma^{-1}}{{(\gamma/\delta)}^{m}\Gamma(\kappa)}\prod_{i=1}^{n}\frac{1}
{\sqrt{2\pi}}S^{-1/2}_{i,p}\]\sqrt{\frac{\pi}{2}}.
$$
This expression holds more generally for $\nu=m+1/2$ or
$\nu=-m-1/2.$
\item[(ii)]If additionally $m=0$, that is $\sum_{j=1}^{p}\theta_{j}=1/2$
and $[\sum_{j=1}^{p}\theta_{j}\lambda_{j}]\Delta=1/2$, then
$$
\Lcr(\X|\vartheta_{p})={\mbox e}^{n\bar{A}\beta} \E\[{\mbox
e}^{-\delta\gamma}\frac{2\gamma^{-1}}{\Gamma(\kappa)}\prod_{i=1}^{n}\frac{1}
{\sqrt{2\pi}}S^{-1/2}_{i,p}\]\sqrt{\frac{\pi}{2}}.
$$
\EndEnumerate \qed
\end{prop}
\Remark Note that superpositioning allows more flexibility in
terms of the parameter values for the constraints
$\sum_{j=1}^{p}\theta_{j}=m+1/2$ and
$[\sum_{j=1}^{p}\theta_{j}\lambda_{j}]\Delta=1/2$. But otherwise
preserves the simplicity of the likelihood as seen in (i) and (ii)
of Proposition 3.7.\EndRemark
\subsection{Randomly sampled times}
From a practical point of view it may be desirable to sample at
uneven or random intervals. See for instance Ait-Sahalia and
Mykland~(2003, 2004). The next result shows that the independence
structure still holds (conditionally) but that the individual
terms are not identically distributed.

\begin{prop}Let $0=\gamma_{0}<\gamma_{1}<\gamma_{2}<\ldots<\gamma_{n}$
denote $n$ random times and define
$\Delta_{i}:=\gamma_{i}-\gamma_{i-1}.$ Define
$\tau_{\theta,i}:=\tau_{\theta}(\gamma_{i})-\tau_{\theta}(\gamma_{i-1})$,
and $r_{i}={\mbox e}^{\lambda \gamma_{i}}$ for $i=1,\ldots, n$,
with $r_{0}=0.$ Then it follows that, conditional on
$(\Delta_{1},\ldots, \Delta_{n})$, for $i=1,\ldots,n,$
$$
\lambda\tau_{\theta,i}=(1-{\mbox e}^{-\lambda \Delta_{i}}){\mbox
e}^{-\lambda
\gamma_{i-1}}[v_{\theta}(0)+\sum_{j=1}^{i-1}r_{j}T_{j}M_{j}]+T_{i}[1-M_{i}]
$$
where $(T_{i},M_{i})$ are conditionally independent pairs
independent of $v_{\theta}(0).$ Additionally, for each fixed $i$,
$T_{i}$ and $M_{i}$ are independent with distributions specified
by $T_{i}\overset{d}=T_{\theta \lambda \Delta_{i}}$ and
$M_{i}\overset{d}=M_{(\theta \lambda \Delta_{i})F_{\lambda
\Delta_{i}}}.$ If the $\Delta_{i}$ for $i=1,\ldots, n$ are
independent then the unconditional distribution of the pairs
$(T_{i},M_{i})$ are independent.
\end{prop}
\subsubsection{Time changed Integrated OU-$\Gamma$ processes}
Notice that the previous proposition places minimal constraints on
the possibly random times $(\gamma_{i}).$ Naturally if one can
easily sample $(\Delta_{1},\ldots, \Delta_{n})$, then this would
lead to models which are amenable to likelihood estimation. These
observations lead us to introduce briefly a  class of time changed
integrated OU processes defined as \Eq
\tau_{\theta}(Z(t))=\int_{0}^{Z(t)}v(s)ds=\lambda^{-1}[(1-{\mbox
e}^{-\lambda Z(t)})v_{\theta}(0)+\int_{0}^{Z(t)}(1-{\mbox
e}^{-\lambda(Z(t)-y)})G_{\theta}(d\lambda y)] \label{TintOU}\EndEq
where $Z$ is any subordinator independent of $G_{\theta}.$ The
next result shows how this model is represented by Proposition
3.8.

\begin{prop}Consider $\tau_{\theta}(Z(t))$ defined as in
~\mref{TintOU}. For $i=1,\ldots,n$, define
$\tau_{\theta,i,Z}:=\tau_{\theta}(Z(i\Delta))-\tau_{\theta}(Z(i-1)\Delta)).$
Then it follows that $\tau_{\theta,i,Z}$ is equivalent to a
specific $\tau_{\theta,i}$ in Proposition 3.8 by setting
$\gamma_{i}=Z(i\Delta)$. Furthermore
$\Delta_{i}=Z(i\Delta)-Z((i-1)\Delta)\overset {d}=Z(\Delta)$ are
iid.\qed
\end{prop}
\Remark The time changed process $\mref{TintOU}$ represents an
extremely rich class of models which adds a great deal of
distributional flexibility to the OU-$\Gamma$ models. As seen from
Proposition 3.9 likelihood analysis for such models is again
easily accomplished. In that case there may be additional unknown
parameters associated with $Z$. For instance selecting $Z$ such
that
$$
\E[{\mbox e}^{-\omega Z(\Delta)}]={\mbox
e}^{-\Delta[(b+\omega)^{1/2}-b^{1/2}]}
$$
Corresponds to the case where $Z(\Delta)$ is an Inverse Gaussian
random variable.\EndRemark \Remark One may also replace $Z(t)$ in
\mref{TintOU} with any tractable increasing process. For instance
one may choose $\tau^{*}_{\alpha}(t)$ to be an integrated
OU-$\Gamma$ process independent of $\tau_{\theta}$\EndRemark
\Remark Leverage type models discussed in BNS pose no extra
difficulties. In the simplest likelihood setting, this translates
into replacing $ X_{i}=\mu\Delta+\beta
\tau_{i}+\tau^{1/2}_{i}\epsilon_{i}$ described in~\mref{BNSx},
with
$$X_{i}=\mu\Delta+ \upsilon T_{i}+\beta
\tau_{i}+\tau^{1/2}_{i}\epsilon_{i}.$$  Where
$T_{i}\overset{d}=T_{\theta a}$ and $\tau_{i}$ is otherwise
related to $T_{i}$ by the representation given in Proposition 3.5.
$\upsilon$ is a real-valued unknown quantity.\EndRemark \Remark We
can extend the OU-$\Gamma$ processes based on the homogeneous
process $G_{\theta}$ to one based on an inhomogeneous Gamma
process $G_{\theta \nu}$, where $\nu$ is an appropriately defined
sigma-finite measure. That is the L\'evy exponent for any positive
function $g$ of $G_{\theta \nu}(g)=\int_{0}^{\infty}g(x)G_{\theta
\nu}(dx)$ is given by $\int_{0}^{\infty}d_{\theta}(\omega
g(x))\nu(dx).$ The volatility process is then defined by
$$v_{\theta \nu}(t)={\mbox e}^{-\lambda t}v(0)+{\mbox e}^{-\lambda
t}\int_{0}^{t}{\mbox e}^{\lambda y}G_{\theta \lambda \nu}(dy)
$$
The process $v_{\theta \nu}(t)$ is stationary only in the
homogeneous case. However the independence properties that we
exploited still hold and one has fairly obvious generalizations of
the results we have presented. An advantage is that this is
another way to increase distributional flexibility.\EndRemark
\subsection{The special nature of the OU-$\Gamma$ process as an SV model}
It is important to note that this independence phenomena,
exhibited in Proposition 3.1,  which allows one to easily describe
the joint structure of $(\tau_{1},\ldots, \tau_{n})$ for a
potential SV model is not only due to the usage of a Gamma process
$G_{\theta}.$ That is to say it will not necessarily be true for
non-OU models based on $G_{\theta}.$ To see this define a moving
average process of the type
$$
\int_{0}^{t}(t-x){\mbox e}^{-(t-x)}G_{\theta}(dx)
$$
It is not difficult to see that the analog of~\mref{pair2} amounts
to  $(\int_{0}^{a}{\mbox e}^{-y}G_{\theta}(dy),\int_{0}^{a}y{\mbox
e}^{-y}G_{\theta}(dy).$ To see the problem first set $H_{a}$ to be
uniform $[0,a]$, and $g_{1}(y)={\mbox e}^{-y}$, and
$g_{2}(y)=y{\mbox e}^{-y}.$ Then it clear that the pair above are
equivalent in distribution to the pair \Eq (T_{\theta a}P_{\theta
aH_{a}}(g_{1}),T_{\theta a}P_{\theta
aH_{a}}(g_{2}))\label{NonOU}\EndEq The good point about this
representation is that the marginal distributional results for
Dirichlet process mean functionals apply. This means, for
instance,  that basically all L\'evy moving average processes that
are driven by a $Z$ which is an FGGC have the property that any
calculation involving a one-dimensional random variable can be
calculated using the marginal distributional results for Dirichlet
process mean functionals. This has an immediate consequence for
option pricing formula based on such models.

However it is quite clear from~~\mref{NonOU} that one can
negotiate the dependence structure in a manner similar to
Proposition 3.1, if and only if $P_{\theta aH_{a}}(g_{1})$ can be
expressed as a function of $P_{\theta aH_{a}}(g_{2})$, which is
not true for this example. This is also why in~\mref{pair2} the
OU-$FGGC$ models we shall discuss do not have the structure
exhibited in Proposition 3.1. In other words $Z(a)$ in that
expression has to have a Gamma distribution. Or more generally
expressible as $T_{\theta a}$ and a function of the other
coordinate.
  Of course the OU-$\Gamma$ is not the only Gamma driven SV model that has the ability to be exactly sampled
 as we did in this section. Another example is the Dykstra and Laud~(1981) type model, see also James~(2005b, p. 1784, eq. (29)),
 which takes the simple form
 $$
 \int_{0}^{t}(t-x)G_{\theta}(dx).
 $$
In this case the analogue of~\mref{pair2} amounts to
$(G_{\theta}(a),\int_{0}^{a}yG_{\theta}(dy)).$

\section{General Likelihoods}
We now proceed to show how one may perform likelihood analysis for
more general $(\tau_{1},\ldots, \tau_{n})$
\subsection{Fourier-Cosine integral representation of the likelihood}
 In order to calculate ~\mref{BNSlik} we use the classical Fourier-Cosine
integral \Eq \frac{1}{\pi}\int_{0}^{\infty}\cos(y|A_{i}|){\mbox
e}^{-\frac{y^{2}\tau_{i}}{2}}dy=\frac{1}{\sqrt{2\pi}}\tau^{-1/2}_{i}{\mbox
e}^{-\frac{A^{2}_{i}}{2\tau_{i}}} \label{key1}.\EndEq This is a
special of the Bessel integral identities known as Weber-Sonine
formula. See for instance Andrews, Askey and Roy~(1999, p.222) and
Watson~(1966, p. 394 eq. (4)) for the identity and also those
references for Bessel functions. It now follows rather immediately
that,

\begin{prop}For the model described by \mref{BNSx}, let
$(\tau_{1},\ldots, \tau_{n})$ have an arbitrary distribution where
the joint Laplace transform has a known form. Then the marginal
likelihood is given by, $$ \Lcr(\X|\vartheta)=\frac{{\mbox
e}^{n{\bar A}\beta }}{\pi^{n}}\int_{{\mathbb R}^{n}_{+}}
\E\[\prod_{i=1}^{n}{\mbox
e}^{-(y_{i}^{2}/2+\beta^{2}/2)\tau_{i}}\]
\prod_{i=1}^{n}\cos(y_{i}|A_{i}|)dy_{i}$$ where
$$\E\[\prod_{i=1}^{n}{\mbox
e}^{-(y_{i}^{2}/2+\beta^{2}/2)\tau_{i}}\]$$ is the joint Laplace
transform of $(\tau_{1},\ldots, \tau_{n})$ evaluated at
$\omega_{i}=y_{i}^{2}/2+\beta^{2}/2$ for $i=1,\ldots,
n.$\qed\end{prop}

The next result which first appears in James~(2005c)[see also
James~(2005b)], which can be thought of an unpublished earlier
version of this manuscript, describes the case where $\tau_{i}$ is
representable as a functional of a Poisson random measure. Since
positive L\'evy processes can be constructed from Poisson random
measures this represents a very rich class.

\begin{prop}Let $N$ denote a Poisson random measure on a Polish space $\Xcr$ with sigma-finite mean intensity
$\nu$, such that for each positive function $g$, the corresponding
random variable $N(g)$ has L\'evy exponent $\Psi(\omega
g)=\int_{\Xcr}(1-{\mbox e}^{-g(x)\omega})\nu(dx).$ Suppose that
$\tau_{i}=N(g_{i})$ for positive functions $(g_{i})$ on $\Xcr$.
Then since
$\sum_{i=1}^{n}N(\omega_{i}g_{i})=N(\sum_{i=1}^{n}\omega_{i}g_{i})$,
it follows that $(\tau_{1},\ldots, \tau_{n})$ has the joint L\'evy
exponent $\Psi(\Omega)=\Psi(\sum_{i=1}^{n}\omega_{i}g_{i}).$ Then
for the model described by \mref{BNSx}, the likelihood is given
by,
$$ \Lcr(\X|\vartheta)=\frac{{\mbox e}^{n{\bar A}\beta
}}{\pi^{n}}\int_{{\mathbb R}^{n}_{+}}{\mbox e}^{-\Psi(\Omega)}
\prod_{i=1}^{n}\cos(y_{i}|A_{i}|)dy_{i}$$ where
$\Omega(x)=\sum_{i=1}^{n}\omega_{i}g_{i}(x)$ with
$\omega_{i}=y_{i}^{2}/2+\beta^{2}/2$ for $i=1,\ldots, n.$\qed
\end{prop}

\Remark Notice that we have stated the result in terms of quite
arbitrary $(\tau_{1},\ldots, \tau_{n}).$ This is because the
expression~\mref{key1} has nothing to do with the distributional
properties of $\tau$.\EndRemark

\Remark Hereafter we set \Eq
\CS(\y|\mu)=\prod_{i=1}^{n}\cos(y_{i}|A_{i}|)
\label{prodcos}\EndEq \EndRemark

\Remark The appearance of integrals involving Bessel functions is
certainly not new to applications in finance as can be seen in the
case of the important work of Yor~(1992) on Asian Options. See
also Carr and Schr\"oder~(2004). \EndRemark

\section{General OU likelihoods} We now apply Proposition 4.1, in the
case of where $(\tau_{1},\ldots,\tau_{n})$ are based on the
integrated OU models described by~\mref{tausimp1}. The task is to
calculate the joint Laplace transform evaluated at
$(\omega_{1},\ldots, \omega_{n}).$ This is straightforward from
the construction given section 2.2. which implies that
$$
\lambda
\sum_{i=1}^{n}\omega_{i}\tau_{i}=s_{1}v(0)+\sum_{l=1}^{n-1}[s_{l+1}\textsc{O}_{l}r_{l}+[Z_{l}-\textsl{O}_{l}]\omega_{l}]+[Z_{n}-\textsl{O}_{n}]\omega_{n}
$$
where for $l=1,\ldots, n$, $ s_{l}=(1-{\mbox e}^{-\lambda
\Delta})[\sum_{i=l}^{n}\omega_{i}{\mbox e}^{-\lambda(i-1)\Delta}].
$ Then it is not difficult to see that the joint Laplace transform
of $(\tau_{1},\ldots, \tau_{n})$ is of the form \Eq
\textsf{L}_{1}(\y|\vartheta)={\mbox e}^{-\varphi(s_{1})}{\mbox
e}^{-\Phi(\omega_{n})}\[\prod_{i=1}^{n-1}{\mbox
e}^{-\Phi(\omega_i|v_{i})}\] \label{jL}\EndEq where terms are
explicitly defined in the next result which gives the likelihood.
\begin{prop}For the model described by \mref{BNSx}, let
$(\tau_{1},\ldots, \tau_{n})$ be defined by the OU models as
in~\mref{tausimp1}. Then the marginal likelihood in~\mref{BNSlik}
is, $$\Lcr(\X|\vartheta)=\frac{{\mbox e}^{n{\bar A}\beta
}}{\pi^{n}}\int_{{\mathbb
R}^{n}_{+}}\textsf{L}_{1}(\y|\vartheta)\CS(\y|\mu)\prod_{i=1}^{n}dy_{i},$$
where $\omega_{i}=y_{i}^{2}/2+\beta^{2}/2$ and
$v_{i}=r_{i}s_{i+1}.$ $\CS(\y|\mu)$ is defined in~\mref{prodcos}
and $\textsc{L}_{1}(\y|\vartheta)$ is the joint Laplace transform
evaluated at $(\omega_{1},\ldots, \omega_{n})$ with form specified
by~\mref{jL}. The L\'evy exponents in~\mref{jL} are specifically
defined as follows, for $a=\lambda \Delta$, \Enumerate
\item[(i)]
$\Phi(\omega_{i}|v_{i})=\int_{{\mbox
e}^{-a}}^{1}\psi(\lambda^{-1}[v_{i}u+\omega_{i}(1-u)])\frac{du}{u},$
for $i=1,\ldots, n-1$
\item[(iii)]
$\Phi(\omega_{n})=\int_{{\mbox
e}^{-a}}^{1}\psi(\lambda^{-1}\omega_{n}(1-u))\frac{du}{u}$
\item[(iv)]
$\varphi(s_{1})=\int_{0}^{1}\psi(s_{1}\lambda^{-1}u)\frac{du}{u},$
is the L\'evy exponent of $v(0)$ evaluated at $s_{1}\lambda^{-1}$.
\EndEnumerate \qed
\end{prop}
Consider now the following result which we will return to in
section 6.
\begin{prop} Consider $\Phi(\omega_{i}|v_{i})$,
$\Phi(\omega_{i})$, and let
$\Lambda(v_{i}|\omega_i)=\Phi(\omega_{i}|v_{i})
-\Phi(\omega_{i}).$ Define
$$D_{\rho}(y,e^{a}y|\omega_{i})=\int_{y}^{\infty}{\mbox
e}^{-\omega_{i}s}\rho(ds)-\int_{y{\mbox e}^{a}}^{\infty}{\mbox
e}^{-\omega_{i}s}\rho(ds)$$ Then \Enumerate
\item[(i)]$\Lambda(v_{i}|\omega_i)=a\int_{0}^{1}\int_{0}^{\infty}(1-{\mbox e}^{-v_{i}us})
{\mbox e}^{-\omega_{i}(1-u)s}\rho(ds)F_{a}(du).$
\item[(ii)]$\Lambda(v_{i}|\omega_i)=\int_{0}^{\infty}(1-{\mbox e}^{-v_{i}y})D_{\rho}(y,e^{a}y|w_{i})y^{-1}{\mbox
e}^{-y}{\mbox e}^{w_{i}y}dy.$
\item[(iii)]It follows that for each fixed $\omega_{i}$,
$\Lambda(t|\omega_{i})$ is the L\'evy exponent, evaluated at $t$,
of an infinitely divisible random variable with L\'evy density
$D_{\rho}(y,e^{a}y|w_{i})y^{-1}{\mbox e}^{-y}{\mbox
e}^{w_{i}y}.$\EndEnumerate
\end{prop}
\Proof Note that $\Lambda(v_{i}|\omega_i)=a\int_{{\mbox
e}^{-a}}^{1}[\psi([v_{i}u+\omega_{i}(1-u)])-\psi([\omega_{i}(1-u)])]
F_{a}(du).$ Statement (i) is simply the L\'evy density
representation of this. Statement (ii) follows by the change of
variable $y=us,$ and exploiting the scale invariance of the
measure $u^{-1}du$. \EndProof This allows one to better understand
the representation of the joint Laplace transform\Eq
\textsf{L}_{1}(\y|\vartheta)={\mbox e}^{-\varphi(s_{1})}{\mbox
e}^{-\Phi(\omega_{n})}\[\prod_{i=1}^{n-1}{\mbox
e}^{-\Phi(\omega_i|v_{i})}\]= {\mbox
e}^{-\varphi(s_{1})}\[\prod_{i=1}^{n}{\mbox
e}^{-\Phi(\omega_{i})}\]\prod_{i=1}^{n-1}{\mbox
e}^{-\Lambda(v_{i}|\omega_{i})} \label{jL1}\EndEq where
$\omega_{i}$ depends only on $y_{i}$ and each $v_{i}$ depends on
$(y_{i+1},\ldots, y_{n}).$ We also note that \Eq
\textsc{L}_{2}(\y|\vartheta):={\mbox
e}^{-\varphi(s_{1})}\[\prod_{i=1}^{n}{\mbox
e}^{-\Phi(\omega_{i})}\] \label{jL2}\EndEq is also a joint Laplace
transform. In fact examining~\mref{jL2} more closely we see that
it is the joint Laplace of a sequence random variables
$(\xi_{1},\ldots,\xi_{n})$, where $ \lambda
\xi_{i}=c_{i}v(0)+[Z_{i}-\textsl{O}_{i}].$ Here $c_{i}=(1-{\mbox
e}^{-a}){\mbox e}^{-a(i-1)}.$ However note that if
$c_{i}=(1-{\mbox e}^{-a})$, then when $v(t)$ is stationary, it
follows that the marginal distribution of this version of
$\xi_{i}$ is equivalent to $\tau_{i}.$ Since we later propose the
use of  joint densities based on~\mref{jL1} and~\mref{jL2} one may
want to vary the value of $c_{i}$ in~\mref{jL2} as this may
increase accuracy.
\section{Some distribution theory for OU-FGGC models}
We have already mentioned that the class of infinitely divisible
random variables which are GGC's are closely linked with Dirichlet
process mean functionals. When the Gamma process has a finite
shape measure say $\theta H$, then every such GGC can be expressed
as $T_{\theta}M_{\theta H}$ where $M_{\theta
H}=\int_{0}^{\infty}xP_{\theta H}(dx)$ is a Dirichlet process. We
will call such GGC's \emph{finite} GGC's or FGGC. One implication
is that one may apply some of the distribution theory we have
developed for the OU-$\Gamma$ to these models.  In this section we
shall assume that $Z$ is derived from a  finite GGC and
demonstrate some nice properties of the corresponding OU process
which are also relevant to sampling likelihoods and option pricing
calculations. First note that if $Z$ is an FGGC then its L\'evy
density and L\'evy exponent are given by \Eq
\psi(\omega)=\int_{0}^{\infty}d_{\theta}(\omega x)H(dx){\mbox {
and }}\theta y^{-1}\int_{0}^{\infty}{\mbox
e}^{-y/r}H(dr)\label{LFGGC}\EndEq where $H$ is a probability
measure

\Remark The term \emph{finite} GGC should not be confused with the
term \emph{finite activity}. That is to say FGGC are
\emph{infinite activity} models as can be seen from their L\'evy
density in ~\mref{LFGGC}. \EndRemark

\subsection{Results for perfect sampling relevant OU-FGGC
components}
\begin{prop} Suppose that $Z$ is a BDLP with
specifications given in ~\mref{LFGGC}. Then the following results
hold. \Enumerate
\item[(i)]In the stationary case the corresponding OU process
$v(t)$ is such that $v(0)$ is a non-finite GGC with L\'evy
exponent
$$
\int_{0}^{\infty}d_{\theta}(\omega
u)S(u)u^{-1}du=\int_{0}^{1}\int_{0}^{\infty}d_{\theta}(\omega
ux)H(dx)u^{-1}du
$$
where $S(u)=\int_{u}^{\infty}H(dy)$ is a survival function.
\item[(ii)]Consider the L\'evy exponent $\Phi(\omega)$ described
in Proposition 4.1. Then in this setting it takes the form \Eq
\Phi(\omega)=\int_{0}^{\infty}d_{\theta a}(\omega r)Q_{a}(dr)
=\int_{{\mbox e}^{-a}}^{1}\int_{0}^{\infty}d_{\theta a}(\omega
(1-u)x)H(dx)F_{a}(du)\label{GGCLevy} \EndEq where $Q_{a}$ is a
probability measure corresponding to the distribution of a random
variable $R=(1-U)W$ where $U$ has distribution $F_{a}$ and $W$ is
independent of $U$ and has distribution $H.$
\item[(iii)]Equivalently $Z(a)-Y_{\theta a}\overset{d}=\int_{0}^{\Delta}(1-{\mbox e}^{-\lambda(\Delta-y)})Z(d\lambda y)$ is a random variable with L\'evy
exponent~\mref{GGCLevy} and hence is a finite GGC and has the
representation $Z(a)-Y_{\theta a}\overset{d}=T_{\theta a}M_{\theta
a Q_{a}}$
\item[(iv)]If the support of $H$ is finite then the support of
$Q_{a}$ is finite.This implies that the perfect simulation method
described by~\mref{upper} and ~\mref{lower} applies to $M_{\theta
aQ_{a}}.$ One may choose $uM_{\theta aQ_{a},-N}$ and $lM_{\theta
aQ_{a},-N}$, according to the upper and lower support points of
$Q_{a}.$ $B_{n,\theta a}$ is \textsc{Beta} $(1,\theta a)$ and
$X_{n}\overset{d}=R=W(1-U)$ has distribution $Q_{a}.$
\item[(v)] For $\theta a=1$ the density of $M_{Q_{a}}$
has the form,
$$
f_{M_{Q_{a}}}(x)=\frac{1}{\pi}\sin(\pi Q_{a}(x)){\mbox
e}^{-\int_{0}^{\infty}\log(|t-x|)Q_{a}(dt)}
$$
where $Q_{a}(x)=\int_{0}^{x}Q_{a}(dt).$ Densities for $\theta a>1$
are obtained by substituting $\theta a=\theta$ and $Q_{a}=H$
in~\mref{Mdensity2}.
\item[(vi)]The Levy exponent of~$Y_{\theta a}\overset{d}=\int_{0}^{\Delta}{\mbox e}^{-\lambda(\Delta-y)}Z(d\lambda y)$ is similar
to ~\mref{GGCLevy} but with $d_{\theta a}(\omega ux)$ in place of
$d_{\theta a}(\omega (1-u)x)$. Hence $Y_{\theta a}\overset
{d}=T_{\theta a} M_{\theta a {\tilde Q}_{a}}$, where ${\tilde
Q}_{a}$ corresponds to the distribution of $\tilde{R}=WU.$
\item[(vi)]$Y_{\theta a}$ converges in distribution to
$v_{\theta}(0)$ as ${\mbox e}^{-a}\rightarrow 0.$\EndEnumerate\qed
\end{prop}
\Proof The proof of (i) and (ii) are obvious by substituting the
form of $\psi$ in~\mref{LFGGC} into~\mref{OLevy} and~\mref{VLevy}.
The remaining results follow as consequences. The density in $(v)$
is obtained from Cifarelli and Regazzini~(1990) or Cifarelli and
Melilli (2001). \EndProof

The next result is rather curious but as we shall show can play a
powerful role in Monte Carlo procedures.
\begin{prop} Consider the setting in Proposition 6.1 then the
L\'evy exponent $\Lambda(t|\omega_{i})$ described in Proposition
5.2 takes the form
$$
\Lambda(t|\omega_{i})=\int_{0}^{\infty}d_{\theta
a}(tr)Q_{a|\omega_{i}}(dr).
$$
where $Q_{a|\omega_i}$ is a probability measure corresponding to a
random variable
$$
R=\frac{UW}{1+W(1-U)\omega_{i}}
$$
where $U$ has distribution $F_{a}$ and $W$ has distribution $H.$
As a consequence results analogous to Proposition 6.1 apply to
this setting.
\end{prop}
\Proof Similar to Proposition 5.2 we examine
$$\Lambda(v_{i}|\omega_i)=a\int_{{\mbox
e}^{-a}}^{1}[\psi([v_{i}u+\omega_{i}(1-u)])-\psi([\omega_{i}(1-u)])]
F_{a}(du).$$ However in this case $
\psi([v_{i}u+\omega_{i}(1-u)])-\psi([\omega_{i}(1-u)]) $ is
equivalent to
$$
\int_{0}^{\infty}[d_{\theta}([v_{i}u+\omega_{i}(1-u)]y)-d_{\theta}([\omega_{i}(1-u)]y)]H(dy)
$$
Now using properties of the natural logarithm it follows that
$$
d_{\theta}([v_{i}u+\omega_{i}(1-u)]y)-d_{\theta}([\omega_{i}(1-u)]y)=d_{\theta}\(\frac{v_{i}uy}{1+\omega_{i}(1-u)y}\)
$$
concluding the result.\EndProof

\subsection{OU-FGGC Monte Carlo Densities}
The next result, whose present importance is that it  can be used
effectively in Monte Carlo simulation procedures, follows
immediately from Propositions 6.1 and 6.2 and standard
augmentation arguments.
\begin{thm} Suppose that the joint Laplace transforms
$\textsc{L}_{1}(\y|\vartheta)$, and hence
$\textsc{L}_{2}(\y|\vartheta)$, satisfies the conditions in
Proposition 6.1 and 6.2. Specify
$\omega_{i}=(y^{2}_{i}+\beta^{2})/2.$ From this, for $i=1,\ldots,
n$, we can let $(T_{i},M_{i})\overset{d}=(T_{\theta a},T_{\theta
a}M_{\theta a Q_{a}})$ denote iid pairs of random variables.
Similarly, independent of the above sequence, define independent
pairs $(G_{i},M_{\omega_{i}})$, where $G_{i}\overset{d}=T_{\theta
a}$ and independent of $G_{i}$, $M_{\omega_i}$ has the
distribution of a mean functional described in Proposition 6.2 for
fixed $\omega_{i}.$ Let
$$\Xi_{1}=(v(0),(T_{i},M_{i}),(G_{i},M_{\omega_{i}}))$$ denote the
joint vector of $4n+1$ independent components, with joint density
$f_{\Xi_{1}}(\cdot|\vartheta, \y).$ Similarly let
$\Xi_{2}=(v(0),(T_{i},M_{i}))$ denote the joint vector of $2n+1$
independent components with density $f_{\Xi_{2}}(\cdot|\vartheta)$
specified by Proposition 6.1 and not depending on $\y.$ Then,
\Enumerate
\item[(i)]
$
\textsc{L}_{1}(\y|\vartheta)=\E[{\mbox
e}^{-s_{1}v(0)}]\prod_{i=1}^{n}\E[{\mbox
e}^{-\omega_{i}T_{i}M_{i}}]\prod_{i=1}^{n}\E[{\mbox
e}^{-v_{i}G_{i}M_{\omega_{i}}}]
$
\item[(ii)]$L_{2}(\y|\vartheta)=\E[{\mbox
e}^{-s_{1}v(0)}]\prod_{i=1}^{n}\E[{\mbox
e}^{-\omega_{i}T_{i}M_{i}}]$
\item[(iii)] Suppose that $\int_{{\mathbb R}^{n}_{+}}
\textsc{L}_{1}(\y|\vartheta) \prod_{i=1}^{n}dy_{i}<\infty$, then
by augmenting the expression in (i) there exists a joint density
of $(\Xi_{1},\Y)$ given by \Eq
f_{\Xi_{1}}(\zeta_{1},\y|\vartheta)\propto {\mbox
e}^{-s_{1}v}\prod_{i=1}^{n}{\mbox
e}^{-\omega_{i}t_{i}m_{i}}\prod_{i=1}^{n}{\mbox
e}^{-v_{i}u_{i}m_{\omega_{i}}}f_{\Xi_{1}}(\zeta_{1}|\vartheta,
\y)\label{L1FGGC}\EndEq where
$\zeta_{1}=(v,(t_{i},m_{i}),(u_i,m_{\omega_{i}})),$ with obvious
meaning.
\item[(iv)]Suppose that $\int_{{\mathbb R}^{n}_{+}}
\textsc{L}_{2}(\y|\vartheta) \prod_{i=1}^{n}dy_{i}<\infty$, then
by augmenting the expression in (ii) there exists a joint density
of $(\Xi_{2},\Y)$ given by $$
f_{\Xi_{2}}(\zeta_{2},\y|\vartheta)\propto {\mbox
e}^{-s_{1}v}\prod_{i=1}^{n}{\mbox
e}^{-\omega_{i}t_{i}m_{i}}f_{\Xi_{2}}(\zeta_{2}|\vartheta, \y),$$
where $\zeta_{2}=(v,(t_{i},m_{i})).$
\item[(v)]Writing $v(0)=T_{\theta}{\tilde M}_{\theta}$ and
integrating out all the Gamma random variables in (iii) it follows
that there exist a joint density of $({\tilde
M}_{\theta},(M_{i}),(M_{\omega_{i}}),\Y)$ given proportional to
\Eq f_{{\tilde
M}_{\theta}}(t){(1+s_{1}t)}^{-\theta}\prod_{i=1}^{n}{(1+\omega_{i}m_{i})}^{-\theta
a}{(1+v_{i}r_{i})}^{-\theta a} f_{M_{i}}(m_{i})
f_{M_{\omega_{i}}}(r_{i}).\label{meanden}\EndEq
 \EndEnumerate
\end{thm}
\subsection{OU-FGGC option pricing densities}
The last result, extends Proposition 3.4 and again is pertinent to
the option pricing formula discussed in BNS(2001a, 6.2) and
Nicolato and Vernados~(2003).

\begin{prop} Let $x^{*}(t)$ be defined by the BDLP $Z$
which is an FGGC with specifications ~\mref{LFGGC}. Additionally,
for $0\leq s<t$, set $\Delta=(t-s)$ and define
$h(\Delta,s)=(1-{\mbox e}^{-\lambda \Delta})v(s)$ and
$\mu^{*}_{s}=\mu\Delta+x^{*}(s)+\beta h(\Delta,s).$ Then the
conditional density of $x^{*}(t)|x^{*}(s),v(s)$ is given by
$$
\int_{0}^{\infty}\phi(x|\mu^{*}_{s}+\beta y,h(\Delta,s)+
y)q_{\theta a}(y)dy
$$
where $q_{\theta a}(y)=\int_{0}^{\infty}\Gcr_{\theta
a}(y|v)f_{M_{\theta aQ_{a}}}(v)dv.$ With the density further
described by the specifications in Proposition 6.1\qed
\end{prop}

\section{Some practical issues for general OU likelihood estimation}
The likelihoods given in Propositions 4.1, 4.2 and 5.1 serve the
purpose of integrating out the infinite-dimensional nuisance
parameters. A natural question is how to exploit these results in
a practical sense. In the forthcoming sections we shall focus on
the OU models but many parts of our discussion can be extended to
more general processes where the joint Laplace transform has an
accessible form. Our goal at minimum will be to discuss ways to
evaluate the likelihood by Monte Carlo procedures. This could then
be used in conjunction with simulated maximum likelihood
estimation or other such techniques. Similar to section 3.4 we
will also be thinking about Bayesian type estimation procedures.
That is, we wish to calculate \Eq
\E\[h(\vartheta)|\X\]=\frac{\int_{\Scr}h(\vartheta)\pi(d\vartheta){{\mbox
e}^{n{\bar A}\beta
}}\textsc{L}_{1}(\y|\vartheta)\CS(\y|\mu)\prod_{i=1}^{n}dy_{i}}{\int_{\Scr}\pi(d\vartheta){{\mbox
e}^{n{\bar A}\beta
}}L_{1}(\y|\vartheta)\CS(\y|\mu)\prod_{i=1}^{n}dy_{i}}
\label{Bayes2}\EndEq where we set $\Scr=({\mathbb
R}^{n}_{+},\Theta).$
\subsection{Calculating L\'evy exponents} In order to utilize
Proposition 5.1 one needs a manageable expression for \Eq
\Phi(\omega_{1}|\omega_{2})= \int_{{\mbox
e}^{-a}}^{1}\psi([\omega_{2}u+\omega_{1}(1-u)])\frac{du}{u}:=a\int_{{\mbox
e}^{-a}}^{1}\psi([\omega_{2}u+\omega_{1}(1-u)])F_{a}(du)\label{jlevy}\EndEq
where $\omega_{1}$ and $\omega_{2}$ just denote two arbitrary
non-negative numbers. We have removed the dependence on the scale
factor $\lambda^{-1},$ which can otherwise be absorbed in
$(\omega_{1},\omega_{2})$. We will assume that $\varphi$ has a
known form. One can see that~\mref{jlevy} is the L\'evy exponent
of the joint distribution of ~\mref{pair2}. However, even if we
wished to try to apply a direct inversion the results for the
OU-$\Gamma$ would suggest that, in general, the joint density of
~\mref{pair2} has a rather non-obvious form. That is to say,
except for the OU-$\Gamma$ case, it is probably just as well to
work directly with~\mref{jlevy}. Now again note importantly that
in order to calculate ~\mref{jlevy} we only need knowledge of
$\psi$ and not the L\'evy density $\rho$ of $Z.$ In many  cases
manual evaluation of~\mref{jlevy} may not be obvious. One can then
resort to numerical methods available in standard mathematical
packages or one can carry out a one time Monte-Carlo approximation
based on the following, somewhat obvious, result.
\begin{prop}Let $U_{l}$ for $l=1,\ldots, B$ denote iid random variables with distribution $F_{a}$. Let
$B\hat{\Phi}(\omega_{1}|\omega_{2})=
\sum_{l=1}^{B}\psi(\omega_{2}U_{l}+\omega_{1}(1-U_{l}))$ Then
$$
\E\[\hat{\Phi}(\omega_{1}|\omega_{2})\]=\Phi(\omega_{1}|\omega_{2})
$$\qed
\end{prop}
Note that our intention is to use a one time calculation of
$\hat{\Phi}(\omega_{1}|\omega_{2})$, based on large $B$,  to get a
highly accurate approximation to $\Phi(\omega_{1}|\omega_{2}).$
Our intent is not to continuously generate different realizations
of $\hat{\Phi}(\omega_{1}|\omega_{2})$ within a loop. In other
words one stores a set of $(U_{l}).$ The remaining sections will
assume that we have been able to get an expression for
$\Phi(\omega_{1}|\omega_{2})$ by some means. \Remark As seen from
our results in section 6 we do not necessarily need to work with
\mref{jlevy} in the case of OU-FGGC models.\EndRemark

\subsection{Monte Carlo method} It is well-known that classical
iid Monte-Carlo, MCMC and SIS  procedures are well-suited to
high-dimensional integrals. However, at first glance, one might
think it is difficult to work with the expressions involving
cosines. Specifically our likelihoods are expressed in terms of $
\CS(\y|\mu)$ which oscillates between positive and negative
values. On the other hand, we note that
$$
|\CS(\y|\mu)|\leq |\cos(y_{1}|A_{1}|)|\leq 1
$$
for all $(y_{1},\ldots, y_{n}),$ which suggests that a product of
cosines is not any more unstable than a single cosine. Monte Carlo
procedures just require a reasonable proposal density and
otherwise deal with terms such  $\CS(\y|\mu)$ in terms of an
expectation $\E[h(\Y)]$ where $h$ depends on $\CS(\y|\mu)$ and
possibly other terms. Accuracy then becomes primarily a function
of the number $B$ of computer iterations. That is, in terms of $B$
Monte Carlo replications. This is in contrast to numerical
techniques which have difficulty handling high dimensions in $n.$
See for instance Liu~(2001), Chen, Shao and Ibrahim~(2000) and
Kong, Liu and Wong~(1997).

The idea of Monte Carlo in the general setting is in principle no
different than that outlined in section 3.4. Except now we will
sample from densities built from $\textsc{L}_{1}(\y|\vartheta)$
and $\textsc{L}_{2}(\y|\vartheta).$  Similar to Theorem 6.1, this
would be possible if its prospective normalizing constant was
finite. Note one can sample these densities without explicit
knowledge of the normalizing constant via MCMC methods. We now
give a description of its normalizing constant.
\begin{prop}Let $\xi_{i}=c_{i}v(0)+[Z_{i}-\textsc{O}_{i}]$,
for $i=1,\dots, n.$ Then
$\max(c_{i}v(0),[Z_{i}-\textsc{O}_{i}])\leq \xi_{i}\leq\tau_{i},$
and the following results hold. \Enumerate
\item[(i)]
$ N_{\vartheta,1}=\int_{{\mathbb R}^{n}_{+}}
\textsc{L}_{1}(\y|\vartheta)
\prod_{i=1}^{n}dy_{i}=\pi^{n}\E\[\prod_{i=1}^{n}\frac{{\mbox
e}^{-\beta^{2}\tau_{i}}}{\sqrt{2\pi \tau_{i}}}\]\leq
\pi^{n}\E\[\prod_{i=1}^{n}\frac{1}{\sqrt{\tau_{i}}}\] $
\item[(ii)]$N_{\vartheta,2}=\int_{{\mathbb R}^{n}_{+}} \textsc{L}_{2}(\y|\vartheta)\prod_{i=1}^{n}
dy_{i}:=\pi^{n}\E\[\prod_{i=1}^{n}\frac{{\mbox
e}^{-\beta^{2}\xi_{i}}}{\sqrt{2\pi \xi_{i}}}\]\leq
\pi^{n}\E\[\prod_{i=1}^{n}\frac{1}{\sqrt{\xi_{i}}}\]$ \qed
\EndEnumerate
\end{prop}
Proposition 7.2 follows by a straightforward argument which can be
seen more clearly in section 8.  We see from Proposition 7.2 that
the prospective normalizing constants are just based on negative
moments of the random variables $\tau_{i}$, $\xi_{i}$, $v(0)$ and
$Z_{i}-\textsc{O}_{i}$ which may or may not exist. One can always
ensure finiteness by adding a small positive constant to any of
the random variables. For instance one uses the model based on
$(\tau_{1}+b,\ldots, \tau_{n}+b)$ for a small $b>0.$ Hereafter we
shall then assume that modification is made if deemed necessary.
This now allows us to describe two possible densities for Monte
Carlo implementation as follows \Eq
{\textsc{N}_{\vartheta,1}}Q_{1}(\y|\vartheta)=
\textsc{L}_{1}(\y|\vartheta) \label{bestden}\EndEq and \Eq
{\textsc{N}_{\vartheta,2}}Q_{2}(\y|\vartheta)=\textsc{L}_{2}(\y|\vartheta).
\label{sbest}\EndEq

Naturally, from the point of view of Monte Carlo (theoretical)
accuracy, $Q_{1}(\y|\vartheta)$ is the most desirable. However,
$Q_{2}(\y|\vartheta)$ is in general easier to sample from. One can
also adjust $Q_{2}(\y|\vartheta)$ further if necessary. Define the
ratio
$$
\Psi(\y|\vartheta)=\frac{Q_{1}(\y|\vartheta)}{Q_{2}(\y|\vartheta)}
=\frac{\textsc{N}_{\vartheta,2}}{\textsc{N}_{\vartheta,1}}\prod_{i=1}^{n}{\mbox
e}^{-\Lambda(v_{i}|\omega_{i})}.
$$

\begin{prop} Consider the densities defined in~\mref{bestden} and~\mref{sbest} and the Bayesian
posterior quantity given in ~\mref{Bayes2}. Additionally let
$\E_{\vartheta,j}$ denote expectation with respect to the
respective joint density of $\Y=(Y_{1},\ldots,Y_{n})$,
$Q_{j}(\y|\vartheta)$ for $j=1,2.$ Define also $\E_{j}$ to denote
expectation with respect to the joint densities
$\pi(\vartheta)Q_{j}(\y|\vartheta)$ for $j=1,2$ Then it follows
that \Enumerate
\item[(i)] $ \Lcr(\X|\vartheta)=N_{\vartheta,1}\frac{{\mbox
e}^{n{\bar A}\beta }}{\pi^{n}}\E_{\vartheta,1}[\CS(\Y|\mu)]$
\item[(ii)]
$ \Lcr(\X|\vartheta)=N_{\vartheta,1}\frac{{\mbox e}^{n{\bar
A}\beta
}}{\pi^{n}}\E_{\vartheta,2}[\CS(\Y|\mu)\Psi(\Y|\vartheta)]$
\item[(iii)]This implies that
\Eq \E\[h(\vartheta)|\X\]=\frac{\E_{1}[h(\vartheta){{\mbox
e}^{n{\bar A}\beta }}N_{\vartheta,1}\CS(\Y|\mu)]}{\E_{1}[{{\mbox
e}^{n{\bar A}\beta
}}N_{\vartheta,1}\CS(\Y|\mu)]}=\frac{\E_{2}[h(\vartheta){{\mbox
e}^{n{\bar A}\beta
}}N_{\vartheta,2}\CS(\Y|\mu)\Psi(\Y|\vartheta)]}{\E_{2}[{{\mbox
e}^{n{\bar A}\beta
}}N_{\vartheta,2}\CS(\Y|\mu)\Psi(\Y|\vartheta)]}
\label{Bayes3}\qed\EndEq
 \EndEnumerate
\end{prop}
Hence a Bayesian approach proceeds similar to section 3.4 by
sampling $(Y_{1,l},\ldots, Y_{n,l},\vartheta_{l})$ for $l=1,\ldots
B$ times from either $\pi(\vartheta)Q_{1}(\y|\vartheta)$ or
$\pi(\vartheta)Q_{1}(\y|\vartheta)$ and put them into appropriate
empirical versions of~\mref{Bayes3}. We now say a few more words
about sampling from the respective densities
\subsubsection{Sampling from $Q_{1}$}
In general an exact expression for the conditional marginals of
say  $Y_{k}|Y_{1},\ldots, Y_{n}$ based on $Q_{1}(\y|\vartheta)$
can be worked out but it is a bit tricky. As such we do not
discuss this. Note that for OU-FGGC models one one can definitely
use Theorem 6.1 to sample from the joint distribution of
$(\Xi_{1},\Y,\vartheta)$ based on ~\mref{L1FGGC} or sampling based
on the density~\mref{meanden}. These methods are facilitated by
the fact that we can use the perfect simulation methods described
in section 2.3.2., with specifications given by Proposition 6.1
and 6.2.
\subsubsection{Sampling from $Q_{2}$}
Notice that in general $Q_{2}(\y|\vartheta)$ has an almost
independent structure and hence a rejection sampling procedure is
straightforward. If however we know the distribution of $v_{0}$ we
can introduce a further augmentation based on
$$
{\mbox e}^{-\varphi(s_{1})}=\int_{0}^{\infty}{\mbox
e}^{-vs_{1}}f_{v(0)}(v)dv.
$$
where again $s_{1}=\sum_{i=1}^{n}c_{i}(y^{2}_{i}+\beta^{2})/2.$
Hence the Monte Carlo procedure can be based on a joint density of
$(\Y,V|\vartheta)|\vartheta$ given as $$Q_{2}(\y,v)\propto
\[\prod_{i=1}^{n}{\mbox e}^{-y^{2}_{i}vc_{i}/2}{\mbox
e}^{-\Phi(\omega_{i})}\]{\mbox
e}^{-v\beta^{2}\sum_{i=1}^{n}c_{i}/2}f_{v(0)}(v).$$ In the OU-FGGC
case we may again use Theorem 6.1 in an obvious way.

\section{General approach} So far we have advocated the idea of
sampling using the joint Laplace transform or some
 variation of that. Since we focused on the BNS models we were able to highlight some nice features. However our claim is that one
can implement similar procedures. This leads us to derive a
similar approach
 that is influenced by some arguments in Devroye~(1986a) but where we do not necessarily sample using the Laplace transform.
That is we give another representation of the likelihood that can
be numerically evaluated via the simulation of random variables.
First let ${\bf p}=(p_{1},\ldots,p_{n})$ denote a vector of
positive numbers and for each $i$, let
$$
H(y_{i}|p_{i})=\frac{2}{\sqrt{2\pi p_{i}}}{\mbox
e}^{-\frac{y^{2}_{i}}{2p_{i}}}{\mbox { for }}y_{i}>0
$$
denote a half Normal density.  Now, notice that $ 0\leq
1-\prod_{i=1}^{n}\cos(y_{i})\leq 2$, and \Eq \int_{{\mathbb
R}^{n}_{+}}\[1-\prod_{i=1}^{n}\cos(y_{i}|A_{i}|)\]
H(y_{i}|p_{i})dy_{i}=1-{\mbox e}^{-\frac{\sum_{i=1}^{n}A^{2}_i
p_{i}}{2}}=C_{n}({\mathbf A},{\mathbf p})\label{cosden} \EndEq
This follows from applications of the Fourier-Cosine identity that
we used in section 4.1. From these facts we describe a joint
density
\begin{prop}Augmenting the expression in~\mref{cosden} leads to a
joint density of an array of positive random variables
$\Y=\{Y_{1,n},\ldots,Y_{n,n}\}$ given by,
$$
r_{n}(\y|{\bf p})=\frac{\[1-\prod_{i=1}^{n}\cos(y_{i}|A
_{i}|)\]\prod_{i=1}^{n}H(y_{i}|p_{i})}{ C_{n}({\mathbf A},{\mathbf
p})}
$$
Equivalently, for $k=1,\ldots, n$, the conditional density of
$Y_{k,n}|Y_{1,n},\ldots, Y_{k-1,n}$ is proportional to
$[1-\lambda_{k}cos(y_{k}|A_{k}|)]H(y_{k}|p_{k})$, where
$\lambda_{k}={\mbox
e}^{-\sum_{i=k+1}^{n}\frac{A^{2}_{i}p_{i}}{2}}\prod_{i=1}^{k-1}\cos(y_{i}|A_{i}|)$
for $k=2,\ldots, n-1$, $\lambda_{1}= {\mbox
e}^{-\sum_{i=2}^{n}\frac{A^{2}_{i}p_{i}}{2}}$, and
$\lambda_{n}=\prod_{i=1}^{n-1}\cos(y_{i}|A_{i}|).$
\end{prop}
Define the function, verified via Fubini's theorem and standard
Normal integration,
$$
\Upsilon_{n}(\vartheta):=\frac{1}{\pi^{n}}\int_{{\mathbb
R}^{n}_{+}}\E\[\prod_{i=1}^{n}{\mbox
e}^{-(y_{i}^{2}/2+\beta^{2}/2)\tau_{i}}\]
\prod_{i=1}^{n}dy_{i}=\E\[\prod_{i=1}^{n}\frac{{\mbox
e}^{-\beta^{2}\tau_{i}}}{\sqrt{2\pi \tau_{i}}}\]\leq
\E\[\prod_{i=1}^{n}\frac{1}{\sqrt{\tau_{i}}}\]
$$

These points lead to following representation of the likelihood.
\begin{prop} Suppose that for fixed $n$,
$\E\[\prod_{i=1}^{n}\frac{1}{\sqrt{\tau_{i}}}\]<\infty$, then the
likelihood in Proposition 4.1 may be written as
$$
{\mbox e}^{{\bar A}\beta}
\[\Upsilon_{n}(\vartheta)-\frac{C_{n}({\mathbf A},{\mathbf
p})}{\pi^{n}}\E\[\Omega(Y_{1,n},\ldots,Y_{n,n}|\vartheta)\]\]
$$
where
$$
\Omega(y_{1},\ldots,y_{n}|\vartheta)=\frac{\E\[\prod_{i=1}^{n}{\mbox
e}^{-(y_{i}^{2}/2+\beta^{2}/2)\tau_{i}}\]}{\prod_{i=1}^{n}H(y_{i}|p_{i})}
$$
and the random vector $\{Y_{1,n},\ldots,Y_{n,n}\}$ has its joint
distribution described by Proposition 8.1.\qed\end{prop}

\Remark Proposition 8.2 shows that one may approximate the
likelihood by simulating random variables described in Proposition
8.1. Such an approach should work well with a Bayesian procedure.
Methods to easily sample the random variables in Proposition 8.1,
may be deduced from Devroye~(1986a, b). In fact, through a
personal communication with Luc Devroye we were informed that one
at time sampling using the conditional distributions in
Proposition 8.1 is routine as it constitutes essentially a
sampling from a Normal density times a factor between $0$ and $2.$
Hence rejection sampling is easy and furthermore the normalizing
factor is not needed. One may also use other densities.\EndRemark

\Remark Note that one needs also to evaluate
$\Upsilon_{n}(\vartheta)$. Of course this can also be done by a
Monte Carlo procedure using the density in Proposition 8.1.
\EndRemark
\section{Examples}
In this section we will present some  examples where we sketch out
a few details related to our exposition. We will not concern
ourselves too much with constants. Note that all the examples
presented are \emph{infinite-activity} processes. In the case
where the distribution of $v(0)$ is not obvious we would simply
approximate it when it is based on OU-FGGC models using
Proposition 6.1, or choose an arbitrary law for $v(0)$ in a more
general setting.
\subsection{OU-Stable}
Suppose that $Z$ is  stable subordinator of index $\alpha$
specified by $\psi(\omega)=\omega^{\alpha}.$ Then it is known, or
otherwise obvious, that $v(t)$ also has a stable law of index
$0<\alpha<1$ with L\'evy exponent
$\omega^{\alpha}\int_{0}^{1}u^{\alpha-1}du.$ Notice that the
L\'evy exponent of the corresponding
$$
\Phi(\omega)=\omega^{a}\int_{{\mbox
e}^{-a}}^{1}{(1-u)}^{\alpha}u^{-1}du
$$
corresponds also to a stable law of index $\alpha$. Here, for
simplicity of presentation, suppressing constants and setting
$\beta=0$ we may use $Q_{2}$ which is based on sampling the joint
Laplace transform \Eq {\mbox
e}^{-{[\sum_{i=1}^{n}y^{2}_{i}]}^{\alpha}}\prod_{i=1}^{n}{\mbox
e}^{-y^{2\alpha}_{i}} \label{StableL}\EndEq Noting the simplicity
of ~\mref{StableL} it is good to recall that in general the
densities of a stable law  are only known in a complicated form.
So here is a case where a Laplace transform approach is perhaps
preferable despite the availability of the relevant densities. A
nice exception to the preceding comment is when $\alpha=1/2$
corresponding to an inverse Gamma law of index $\alpha=1/2$.
However in that case ~\mref{StableL} is $${\mbox
e}^{-{[\sum_{i=1}^{n}y^{2}_{i}]}^{1/2}}\prod_{i=1}^{n}{\mbox
e}^{-y_{i}}.$$ For further simplification we may use the
augmentation procedure described in section 7.2.2 applied
to~\mref{StableL} to get $$\prod_{i=1}^{n}{\mbox
e}^{-y^{2\alpha}_{i}}{\mbox e}^{-vy^{2}_{i}}f_{\alpha}(v)$$ where
$f_{\alpha}$ corresponds to a stable density. Note that although
the stable density can be complicated there are many routines
available to easily sample stable random variables. \Remark The
Stable law process produces a log price process with heavy tails
which may not be desirable for all applications. However see the
work of Carr and Wu~(2003). Additionally, we note that it would
not be tremendously difficult to use $Q_{1}$ in this case.
\EndRemark
\subsection{IG-OU}This example is based on the calculations given in Barndorff-Nielsen and Shephard~(2003, p. 292) where
$v(t)$ has an Inverse Gaussian distribution. Here, letting
$C_{1},C_{2}$ denote constants and setting $\beta =0$, by
BNS(2003, eq. (54)) one has
$$
\Phi(\omega)=-y^{2}_{i}C_{1}\int_{0}^{1-{\mbox
e}^{-a}}{(1-u)}^{-1}u{(1+C_{2}y^{2}_{i}u)}^{-1/2}du.
$$
BNS(2003) show that this can  be written in terms of the
hyperbolic arc-tangent function[see also Nicolato and
Vernardos~(2001) and Carr, Geman, Madan and Yor~(2003)], we do not
repeat that here. Note however by using the fact that $v(t)$ has a
Inverse Gaussian distribution one can work with the augmented
version of $Q_{2}$ which is proportional to
$$
v^{-3/2}{\mbox e}^{-\frac{1}{2}[\gamma^{2} v+\delta^{2}
v^{-1}]}\prod_{i=1}^{n}{\mbox e}^{-y^{2}_i v}{\mbox
e}^{-\Phi(\omega_{i})}
$$
for appropriate values of $\gamma$ and $\delta$ and is not
difficult to sample from.
\subsection{OU-LogNormal}
Suppose that $Z$ is based on a LogNormal distribution with density
$$
f(x)=\frac{1}{\sqrt{2\pi}}\frac{1}{x}{\mbox
e}^{-\frac{1}{2}(\log(x))^{2}}{\mbox { for }}x>0.
$$
We have chosen this example because, despite the fact that it has
a density with a nice closed form, its corresponding L\'evy
density $\rho$ is unknown. Despite this we can still use a sampler
based on $Q_{2}$. This is because its L\'evy exponent is given by
$$
\psi(\omega)=-\log \[\int_{0}^{\infty}{\mbox e}^{-\omega
x}\frac{1}{\sqrt{2\pi}}\frac{1}{x}{\mbox
e}^{-\frac{1}{2}(\log(x))^{2}}dx\].
$$
This can be numerically approximated hence the relevant quantities
$\Phi(\omega)$ can then be numerically approximated. Again this
approximation should be done before the main Monte-Carlo procedure
is used. Note that one would find it difficult or impossible to
employ a series approximation in this case, as it depends on
knowledge of $\rho.$ \subsection{OU-FGGC where $H$ is the Arcsine
distribution}Here we close with one of the more interesting
examples of known FGGC models. In this setting let $H$ be the
Arcsine law, that is there is a corresponding random variable $W$
which is $\textsc{Beta}(1/2,1/2)$. Cifarelli and Melilli~(2000)
show that in this setting for all $\theta>0$, $M_{\theta H}$ is
$\textsc{Beta}(\theta+1/2,\theta+1/2)$. Hence
$Z(t)\overset{d}=T_{\theta t}B_{(\theta t+1/2,\theta t+1/2)},$
where here $B_{(\theta t+1/2,\theta t+1/2)}$ means a beta random
variable with parameters indicated in the subscript.  In this case
the distribution of $R=(1-U)W$ described in Proposition 6.1 has
bounded support on $[0,1]$. Hence we may apply the perfect sampler
both for option pricing and Monte Carlo methods specified
according to Theorem 6.1 and Propositions 6.1 and 6.2. That is
apply section 2.3.2 to sample from $(\Xi_{1},\Y)$. To be clear
given $Y_{i}$ one may draw $M_{\omega_{i}}$ from
$f_{M_{\omega_{i}}}$ by using Proposition 6.2 and creating
$uM_{\theta a,-N}=1$ and $lM_{\theta a,w_{i},-N}=0$, $B_{n,\theta
a}$ is \textsc{Beta} $(1,\theta a)$ and
$X_{n}\overset{d}=UW/(1+W(1-U)w_{i}).$ Then one draws $W$ from the
Arscine law and $U$ from $F_{a}$ to get $X_{(n)}$. Draws from more
complex densities for $M_{w_{i}}$ can then be obtained by other
standard methods. One can also work with the exact form of the
densities via Cifarelli and Regazzini~(1990).

\vskip0.2in \centerline{\Heading References} \vskip0.2in \tenrm
\def\smc{\tensmc}
\def\sl{\tensl}
\def\bf{\tenbold}
\baselineskip0.15in

\Ref \by A\"it-Sahalia, Y., Mykland, P. A. \yr 2003 \paper The
effects of random and discrete sampling when estimating
continuous-time diffusions \jour Econometrica \vol 71 \pages
483-549\EndRef

 \Ref \by A\"it-Sahalia, Y., Mykland, P. A. \yr
2004 \paper Estimators of diffusions with randomly spaced discrete
observations: a general theory \jour \AnnStat \vol 32 \pages
2186-2222\EndRef

\Ref \by Andrews, G., Askey, R. and  Roy, R. \yr 1999 \book
Special functions. Encyclopedia of Mathematics and its
Applications, 71 \publ Cambridge University Press \publaddr
Cambridge \EndRef

\Ref \by Barndorff-Nielsen, O.E. and Shephard, N. \yr 2001a \paper
Ornstein-Uhlenbeck-based models and some of their uses in
financial economics \jour \JRSSB \vol 63 \pages 167-241 \EndRef
\Ref \by Barndorff-Nielsen, O.E. and Shephard, N. \yr 2001b \paper
Modelling by L\'evy processes for financial econometrics. In
L\'evy processes. Theory and applications. Edited by Ole E.
Barndorff-Nielsen, Thomas Mikosch and Sidney I. Resnick. p.
283-318. Birkh\"auser Boston, Inc., Boston, MA \EndRef

\Ref \by Barndorff-Nielsen, O. E. and Shephard, N.\yr 2003 \paper
Integrated OU processes and non-Gaussian OU-based stochastic
volatility models \jour Scand. J. Statist. \vol 30 \pages 277-295
\EndRef

\Ref \by Benth, F. E., Karlsen, K. H. and Reikvam, K. \yr 2003
\paper Merton's portfolio optimization problem in a Black and
Scholes market with non-Gaussian stochastic volatility of
Ornstein-Uhlenbeck type \jour Math. Finance \vol 13 \pages 215-244
\EndRef

\Ref \by Black, F. and Scholes, M. \yr 1973 \paper The pricing of
options and corporate liabilities \jour J. Polit. Econ. \vol 81
\pages 637-654 \EndRef

\Ref \by Bondesson, L. \yr 1979 \paper A general result on
infinite divisibility \jour \AnnProb \vol 7 \pages 965-979 \EndRef

\Ref \by Bondesson, L. \yr 1992 \paper Generalized gamma
convolutions and related classes of distributions and densities.
Lecture Notes in Statistics, 76. Springer-Verlag, New York \EndRef

\Ref \by Carr, P., Geman, H., Madan, D.B. and Yor, M. \yr 2003
\paper Stochastic volatility for L\'evy processes \jour Math.
Finance \vol 13 \pages 345-382 \EndRef

\Ref \by Carr, P., Geman, H., Madan, D.B. and Yor, M. \yr 2005
\paper Self-Decomposability and Option Pricing. Math. Finance to
appear\EndRef

\Ref \by Carr, P. and Schr\"oder, M. \yr 2004 \paper Bessel
processes, the integral of geometric Brownian motion, and Asian
options \jour Theor. Probab. Appl. \vol 48 \pages 400-425 \EndRef

\Ref \by Carr, P. and Wu, L. \yr 2003 \paper The finite moment log
stable process and option pricing \jour Journal of Finance \vol 58
\pages 753-778\EndRef

\Ref \by Carr, P. and Wu, L. \yr 2004 \paper Time-changed L\'evy
processes and option pricing \jour Journal of Financial Economics
\vol 71 \pages 113-141\EndRef

\Ref \by Chen, M-H., Shao, Q-M.,  Ibrahim, J.G. \yr 2000 \book
Monte Carlo methods in Bayesian computation. Springer Series in
Statistics.\publ Springer-Verlag \publaddr New York\EndRef

\Ref \by    Cifarelli, D. M. and Melilli, E. \yr    2000 \paper
Some new results for Dirichlet priors \jour \AnnStat \vol 28
\pages 1390-1413\EndRef

\Ref \by Cifarelli, D. M. and Regazzini, E.\yr 1990 \paper
Distribution functions of means of a Dirichlet process \jour
\AnnStat \vol 18 \pages 429-442\EndRef

\Ref \by Devroye, L. \yr 1986a \paper An automatic method for
generating random variates with a given characteristic function
\jour SIAM J. Appl. Math. \vol 46 \pages 698-719\EndRef

\Ref \by Devroye, L. \yr 1986b \book Nonuniform random variate
generation. \publ Springer-Verlag \publaddr New York\EndRef \Ref

\by Diaconis, P. and Freedman, D. A. \yr 1999 \paper Iterated
random functions \jour Siam Rev. \vol 41 \pages 45-76 \EndRef

\Ref \by    Diaconis, P. and Kemperman, J. \yr    1996 \paper Some
new tools for Dirichlet priors. Bayesian Statistics 5 (J.M.
Bernardo, J.O. Berger, A.P. Dawid and A.F.M. Smith eds.), Oxford
University Press, pp. 97-106 \EndRef

\Ref \by Duan,  J.\yr 1995 \paper The GARCH option pricing model
\jour Math. Finance \vol 5 \pages 13-32 \EndRef

\Ref \by Duffie, D., Pan, J. and Singleton, K.,\yr 2000 \paper
Transform Analysis and Asset Pricing for Affine Jump Diffusions
\jour Econometrica \vol 68 \pages 1343-1376 \EndRef

\Ref \by Dykstra, R. L. and Laud, P. W. \yr 1981 \paper A Bayesian
nonparametric approach to reliability \jour \AnnStat \vol 9 \pages
356-367 \EndRef

\Ref \by Eberlein, E. \yr 2001 \paper Application of generalized
hyperbolic L\'evy motions to finance. In L\'evy processes. Theory
and applications. Edited by Ole E. Barndorff-Nielsen, Thomas
Mikosch and Sidney I. Resnick. p. 319-336. Birkh\"auser Boston,
Inc., Boston, MA \EndRef

\Ref \by Engle, R. F.\yr 1982 \paper Autoregressive conditional
heteroscedasticity with estimates of the variance of United
Kingdom inflation \jour Econometrica \vol 50 \pages 987-1007
\EndRef

\Ref \by Eraker, B., Johannes, M., and Polson, N. \yr 2003 \paper
The impact of jumps in volatility and returns. \jour Journal of
Finance \vol 68 \pages 1269-1300 \EndRef

\Ref \by Flajolet, P. and Sedgewick, R. \yr 2006 \paper Analytic
Combinatorics. Book to appear. Chapters available at\\
http://algo.inria.fr/flajolet/Publications/books.html\EndRef

 \Ref \by Guglielmi, A., Holmes, C.C., Walker, S.G. \paper
Perfect simulation involving functionals of a Dirichlet process
\jour  J. Comput. Graph. Statist. \vol 11 \pages 306-310 \EndRef

\Ref \by Griffin, J. and Steel, M. \yr 2005 \paper Stochastic
Volatility Inference with non-Gaussian Ornstein-Uhlenbeck
Processes for Stochastic Volatility forthcoming Journal of
Econometrics\EndRef

\Ref \by Guglielmi, A., Holmes, C.C., Walker, S.G. \paper Perfect
simulation involving functionals of a Dirichlet process \jour  J.
Comput. Graph. Statist. \vol 11 \pages 306-310 \EndRef

\Ref \by Hjort, N. L., and Ongaro, A. \yr 2005 \paper Exact
inference for random Dirichlet means. \jour Stat. Inference Stoch.
Process. \vol 8 \pages 227-254 \EndRef

\Ref \by James, L.F.  \yr 2005a \paper Functionals of Dirichlet
processes, the Cifarelli-Regazzini identity and Beta-Gamma
processes \jour \AnnStat \vol 33 \pages 647-660\EndRef

\Ref \by James, L.F. \yr 2005b \paper Bayesian Poisson process
partition calculus with an application to Bayesian L\'evy moving
averages. \jour \AnnStat \vol 33 \pages 1771-1799\EndRef

\Ref \by James, L.F. \yr 2005c \paper Analysis of a class of
likelihood based continuous time stochastic volatility models
including Ornstein-Uhlenbeck models in financial economics.
arXiv:math.ST/0503055 \EndRef

\Ref \by Jeanblanc, M., Pitman, J. and Yor, M. \yr 2002 \paper
Self-similar processes with independent increments associated with
L\'evy and Bessel \jour Stochastic Process. Appl. \vol 100 \pages
223-231\EndRef

\Ref \by Jurek, Z.J., Vervaat, W. \yr 1983 An integral
representation for self-decomposable Banach space valued random
variables \jour \jour Z. Wahrsch. Verw. Gebiete \vol 62 \pages
247-262\EndRef

\Ref \by Kong, A.,  Liu, J.S., and  Wong, W. H. \yr 1997 \paper
The properties of the cross-match estimate and split sampling
\jour \AnnStat \vol 25 \pages 2410-2432\EndRef

\Ref \by Liu, J.S. \yr 2001 \book Monte Carlo strategies in
scientific computing. Springer Series in Statistics \publ
Springer-Verlag \publaddr New York\EndRef

\Ref \by Lukacs, E.A. \yr 1955 \paper A characterization of the
gamma distribution \jour Ann. Math. Statist. \vol 26 \pages
319-324\EndRef

\Ref \by Madan, D., Carr, P. and Chang, E. \yr 1998 \paper The
variance gamma process and option pricing \jour European Finance
Rev. \vol 2 \pages 79-105\EndRef

\Ref \by Merton, R.~C. \yr 1973 \paper Theory of rational option
pricing \jour Bell J. Econ. Mgemt. Sci. \vol 4 \pages 141-183
\EndRef

\Ref \by Maximon, L. C. \yr 2003 \paper The dilogarithm function
for complex argument \jour R. Soc. Lond. Proc. Ser. A Math. Phys.
Eng. Sci. \vol 459 \pages 2807--2819\EndRef

\Ref \by Nicolato, E. and Venardos, E. \yr 2003 \paper Option
pricing in stochastic volatility models of the Ornstein-Uhlenbeck
type \jour Math. Finance \vol 13 \pages 445-466 \EndRef

\Ref \by  Pitman, J., \yr 1999 \paper Brownian motion, bridge,
excursion, and meander characterized by sampling at independent
uniform times\jour Electron. J. Probab. \vol 4 \pages 1-33\EndRef

\Ref \by Propp, J.G. and Wilson, D. B. \yr 1996 \paper Exact
sampling with coupled Markov chains and applications to
statistical mechanics \jour Random Structures Algorithms \vol 9
\pages 223-252 \EndRef

\Ref \by Roberts, G. O., Papaspiliopoulos, O. and Dellaportas, P.
\yr 2004 \paper Bayesian inference for non-Gaussian
Ornstein-Uhlenbeck stochastic volatility processes \JRSSB \vol 66
\pages 369-393 \EndRef

\Ref \by Sato, K. \yr 1999 \book L\'evy processes and infinitely
divisible distributions. Translated from the 1990 Japanese
original. Cambridge Studies in Advanced Mathematics, 68 \publ
Cambridge University Press \publaddr Cambridge \EndRef

\Ref \by Thorin, O.\yr 1977 \paper On the infinite divisibility of
the lognormal distribution \jour Scand. Actuar. J. \vol 3 \pages
121-148 \EndRef

\Ref \by Watson, G. N. \yr 1966 \book A treatise on the theory of
Bessel functions. Paperback Edition. Cambridge Mathematical
Library. Cambridge University Press, Cambridge \EndRef

\Ref \by Wolfe, S. J. \yr 1982 \paper On a continuous analogue of
the stochastic difference equation $X\sb{n}=\rho
X\sb{n-1}+B\sb{n}$ \jour Stochastic Process. Appl. \vol 12 \pages
301-312\EndRef

\Ref \by Vershik, A.M., Yor, M. and Tsilevich, N.V.  \yr 2004
\paper On the Markov-Krein identity and quasi-invariance of the
gamma process \jour J. Math. Sci. \vol   121 \pages 2303-2310
\EndRef

\Ref \by Yor, M. \yr 1992 \paper On some exponential functionals
of Brownian motion \jour Adv. in Appl. Probab. \vol 24 \pages
509-531 \EndRef
\medskip
\smc

\Tabular{ll}

Lancelot F. James\\
The Hong Kong University of Science and Technology\\
Department of Information and Systems Management\\
Clear Water Bay, Kowloon\\
Hong Kong\\
\rm lancelot\at ust.hk\\
\EndTabular

\end{document}